\documentclass[smallextended]{svjour3}

\usepackage{graphicx}
\usepackage{amssymb,amsmath,amsopn,amstext,amsfonts,color,fullpage,hyperref,xcolor}

\usepackage{graphicx,pgf,tikz}
\usetikzlibrary{fadings,shapes,shapes.misc}

\usepackage[boxed,commentsnumbered]{algorithm2e}
\usepackage{rotating}

\makeatletter
\newcommand\footnoteref[1]{\protected@xdef\@thefnmark{\ref{#1}}\@footnotemark}
\makeatother

\def\R{{\mathbb{R}}}

\newcommand{\tm}{\tilde{m}}
\newcommand{\tn}{\tilde{n}}

\newcommand{\real}{{\mathbb{R}}}
\newcommand{\bin}{{\mathbb{B}}}

\newcommand{\one}{\mathbf{1}}

\begin{document}

\bibliographystyle{spmpsci}

\title{Solving Linear Programs with Complementarity Constraints using
Branch-and-Cut\thanks{Yu and Mitchell were
supported in part by the Air Force Office of Sponsored Research under grants FA9550-08-1-0081 and FA9550-11-1-0260
and by the National Science Foundation under Grant Numbers CMMI-1334327
and DMS-1736326.
Pang was
supported in part by the National Science Foundation
			under Grant Number CMMI-1333902
			and by the Air Force Office of Scientific Research
			under Grant Number FA9550-11-1-0151.
}
}

\author{Bin Yu        \and
        John E. Mitchell   \and  Jong-Shi Pang 
}

\institute{B. Yu \at
              BNSF Railway,
Fort Worth, TX     \\   
           \and
           J.E. Mitchell \at
          Department of Mathematical Sciences,
    Rensselaer Polytechnic Institute,
    Troy,  NY 12180.
    \email{mitchj@rpi.edu}
    \url{http://www.rpi.edu/~mitchj}   \\
    \and
    J.S. Pang \at
    Department of Industrial and Systems Engineering,
    University of Southern  California,
     Los Angeles, CA 90089, USA.
    \email{jongship@usc.edu}
}

\date{Received: date / Accepted: date}

\maketitle

\begin{abstract}
A linear program with linear complementarity constraints (LPCC)
requires the minimization of a linear objective over a set 
of linear constraints together with additional linear complementarity constraints.
This class has emerged as a modeling paradigm for a broad
collection of problems, including bilevel programs,
Stackelberg games, inverse quadratic programs,
and problems involving equilibrium constraints.
The presence of the complementarity constraints
results in a nonconvex optimization problem.
We develop a branch-and-cut algorithm to find a global optimum for
this class of optimization problems, where
we branch directly on complementarities.
We develop branching rules and feasibility recovery procedures
and demonstrate their computational effectiveness
in a comparison with CPLEX.
The implementation builds on CPLEX through the use of callback routines.
The computational results show that our approach is a strong alternative
to constructing an integer programming formulation using big-$M$ terms
to represent bounds for variables,
with testing conducted on general LPCCs as well as on instances
generated from bilevel programs with convex quadratic lower level problems.
\end{abstract}

\section{Introduction\label{sec:introduction}}

A linear program with linear complementarity constraints (LPCC), which minimizes a linear objective function over a set of linear constraints with additional linear complementarity constraints, is a non-convex, disjunctive optimization problem.
In \S\ref{sec: Statement of the Problem}, we present the mathematical formulation of the general LPCC we use throughout this paper. 
In \S\ref{sec: Previous Work for Solving LPCCs}, various existing algorithms designed for solving LPCCs are reviewed. Most of these existing methods are only able to obtain a stationary solution and incapable of ascertaining the quality of the solution. This is the major drawback for the existing solvers. In this paper, we mainly focus on finding the global resolution of the LPCC, and we achieve this goal through two steps:
\begin{description}
\item[Step 1]Study various valid constraints by exploiting the complementarity constraints directly, and evaluate the benefit of these constraints on the value of the linear relaxation of the LPCC.
We have previously discussed valid constraints for the LPCC in \cite{JEMPYuMOPTA},
and we briefly recap these constraints in~\S\ref{sec:  Cutting Planes Generation and Selection}.
\item[Step 2]Propose a branch-and-cut algorithm to globally solve the LPCC problem, where cuts are derived from various valid constraints studied in Step 1 and branching is imposed on the complementarity constraints. A general LPCC solver has been developed based on this branch-and-cut approach, and it is able to compete with the existing MIP-based solvers like CPLEX.
\end{description}
The branch-and-cut algorithm is introduced in
\S\ref{sec:LPCC using branch-and-cut}, where we also outline the rest of the paper.

\subsection{Statement of the Problem\label{sec: Statement of the Problem}}

We consider a general formulation of the LPCC in the form suggested by Pang and Fukushima~\cite{pangfukushima1}. 
Given vectors and matrices:
$c \in \R^n$, $d \in \R^m$, $b \in \R^k$, $q \in \R^m$,
$A \in \R^{k \times n}$, $B \in \R^{k \times m}$, $N \in \R^{m \times n}$ and $M \in \R^{m \times m}$,
the LPCC is to find $(x,y,w) \in \R^n \times \R^m \times \R^m$ in order to solve to global optimality
\begin{equation} \label{eq:general LPCC}
\begin{array}{ll}
\displaystyle{
{\operatornamewithlimits{\mbox{minimize}}_{(x,y,w)}}
} & c^Tx + d^Ty\\ 
\mbox{subject to} & Ax + By  \geq  b \\
& x \, \geq \, 0  \\
\mbox{and} & 0  \leq  y  \perp w:= q + Nx + My  \geq  0
\end{array} \end{equation}
where $a \perp b$ denotes perpendicularity between vectors $a$ and $b$, i.e., $a^Tb=0$.
Without the orthogonality condition $y  \perp  w$, the LPCC is a linear program (LP).
The global resolution of the LPCC means the generation of a certificate
showing that the problem is in one of its 3 possible states: (a) it is infeasible,
(b) it is feasible but unbounded below, or (c) it attains a finite optimal solution.
Note that problem~(\ref{eq:general LPCC}) is equivalent to $2^m$ linear programs obtained by making
each possible assignment for the complementary variables: either $y_i=0$ or $w_i=0$ for each $i=1,\ldots,m$;
hence, it is not possible for an LPCC to have a finite optimal value that is not attained.

If the feasible regions for $y$ and $w$ are bounded then there exist diagonal matrices
$\Theta^y$ and $\Theta^w$ with diagonal entries $\theta^y_i$ and $\theta^w_i$
and problem (\ref{eq:general LPCC}) can be formulated as a mixed integer program:
\begin{equation} \label{eq:MIP LPCC}
\begin{array}{ll}
\displaystyle{
{\operatornamewithlimits{\mbox{minimize}}_{(x,y,z)}}
} & c^Tx + d^Ty\\ 
\mbox{subject to} & Ax + By  \geq  b \\ 
 & x \, \geq \, 0  \\
& 0  \leq  y \leq \Theta^y z \\
 & 0 \leq  q + Nx + My  \leq  \Theta^w(\textbf{1}-z)\\
\mbox{and} & z\in\{ 0,1\}^m
\end{array} \end{equation}
The obvious drawback of this formulation is that in order to find $\theta^y_i$ and $\theta^w_i$
we need to compute valid upper bounds of $y_i$ and $w_i$, not to mention such upper bounds may not exist
if the feasible regions for $y$ and/or $w$ are unbounded.
To avoid this drawback, in this paper, we present a branch-and-cut algorithm which branches on the complementarity constraint directly.
Previous work on branching on complementarity constraints
includes~\cite{bealetomlin1970,farias2001,hooker_osorio}.

Problem (\ref{eq:general LPCC}) generalizes the standard linear complementarity problem (LCP)~\cite{CPStone92}:  $0 \leq y \perp q + My \geq 0$,
so the LPCC is NP-Hard. Moreover, 
affine variational constraints also lead to the problem (\ref{eq:general LPCC})~\cite{luo7}.
Applications of the LPCC are surveyed in~\cite{jinghu3}.
Among these applications, complementarity constraints play three principal roles during the modelling process: 
\begin{enumerate}
\item Modelling KKT optimality conditions that must be satisfied by some of the variables. Such applications include hierarchical optimization such as Stackelberg games~\cite{stackelberg1},
inverse convex quadratic programs,
indefinite quadratic programs~\cite{GTomasin73,HMPang08},
and cross-validated support vector regression~\cite{KPBennett2011,yuchingSVM}.
\item Modelling equilibrium constraints. See for example the texts~\cite{facchinei2,luo7},
the survey article~\cite{Pang2010},
or a recent paper on market equilibrium in electric power markets~\cite{MLCPanjos}.
\item Modelling certain logical conditions that are required by some practical optimization problems. Such applications include non-convex piecewise linear optimization,
quantile minimization~\cite{pang2}, and $\ell_0$-minimization~\cite{bkSchwartz2014,FMPSWachter13}.
\end{enumerate}

\subsection{Previous Work on Solving LPCCs\label{sec: Previous Work for Solving LPCCs}}

Research on algorithms for solving an LPCC can be divided into two main areas:
one concerns the development of globally convergent algorithms with a guarantee of finding a suitable stationary point;
the other concerns the development of exact algorithms for global resolution of an LPCC.
See the survey~\cite{judiceTOP2012} for a more detailed review.

 It is noted that the methods which are able to solve a general LCP can also be extended to solve an LPCC by using the so called sequential LCP Method. Such a procedure can be found in detail
in~\cite{judicefaustino1}. 
 A complementary pivoting algorithm for an LPCC is an extension of a pivoting algorithm for LCP which handles linear complementarity constraints just as the classic simplex algorithm for linear programs. Such algorithms usually perform in this way: start from a feasible solution, maintain feasibility for all iterations and try to improve the objective function. Under certain constraint qualifications, these methods guarantee convergence to a certain stationary solution. The references~\cite{fangleyffermunson1,fukushima1,izmailovsolodov09}
study and implement this type of method to solve the general LPCC. Another way to get a stationary point is through a so called regularization framework~\cite{scholtes2}: construct a sequence of relaxed problems controlled by some parameter, then obtain a sequence of solutions which converge to a stationary point when the parameter goes to the limit.  Each regularized relaxed problem is solved by 
an NLP based algorithm such as an interior point method. One method of regularization is to introduce a positive parameter $\phi$ and relax the complementarity constraints in problem (\ref{eq:general LPCC}) using either  $\{y,w\geq 0, y^Tw\leq \phi\}$ or some other approach~\cite{bkSchwartz2014,FMPSWachter13}.
An alternative is to put a penalty for violation of the complementarity constraints into the objective, and gradually update the penalty to infinity~\cite{leyffer4}.
A homotopy method has also been proposed~\cite{watsonBME2013}.
The obvious drawback of these methods is that they are incapable of ascertaining the quality of the computed solution. 

The methods for  global resolution of an LPCC are mainly based on an enumerative scheme.   
Several branch-and-bound methods have been proposed for solving an LPCC derived from a bilevel linear program.
Bard and Moore~\cite{bardmoore1} proposed a pure branch-and-bound method for solving bilevel linear programs. Hansen et al.~\cite{hansenjaumardsavard1}  enhance this branch-and-bound scheme by exploiting the necessary optimality condition of the inner problem. As opposed to a branch-and-bound method, the references~\cite{ibaraki2,jeroslow1} study alternative ways to solve an LPCC by using a cutting plane method. 
Audet et al.~\cite{audet6} proposed a branch-and-cut algorithm for solving bilevel linear programs.
An RLT method for finding a feasible solution to a problem with both binary and complementarity
constraints is proposed in~\cite{MLCPanjos}.
It follows from the results of \cite{lijie2} that an LPCC can be lifted to an equivalent
convex optimization problem so it can in principle be solved globally using a convex optimization algorithm;
the drawback to this approach is that the convexity is over the cone of completely positive matrices
which is hard to work with computationally.

It is noted that most of the existing methods for global resolution of the LPCC
presume the LPCC has a finite optimal value, and this limitation was not resolved until the
paper~\cite{jinghu1}. In that paper, the authors proposed a minimax integer programming formulation of the LPCC, and solve this system using a Benders decomposition method.
The method was extended to quadratic programs with complementarity
constraints in~\cite{lijie1}.
A branching scheme for determining boundedness of the optimal
value of a linear program with a bilinear objective function was
proposed in~\cite{AHJSavard99}.

The success of the Benders decomposition method
\cite{lijie1,jinghu1}
heavily depends on a so called sparsification process.
If the sparsification process is not successful,
in the worst case it will be necessary to check every piece of the LPCC.
In this paper, we alternatively use a specialized
branch-and-cut scheme which is a more systematic enumerative process to get the global resolution of the LPCC,
and our algorithm is also able to characterize infeasible and unbounded LPCC problems as well as solve problems with finite optimal value.
Moreover we also discuss various valid constraints for the LPCC by exploiting the complementarity structure;
this topic has not been fully exploited in the literature for studying the LPCC.

The complementarity structure of an LPCC can be generalized to SOS1 constraints,
a type of special ordered set constraint requiring that at most one of a set of variables
is nonzero.
Recent work on branch-and-cut approaches to problems with SOS1
constraints include~\cite{farias8,fischerpfetsch2015}.
De Farias et al.~\cite{farias8} considered problems where all the coefficients are
nonnegative and their emphasis is on possible families of cutting planes
using a sequential lifting procedure.
Fischer and Pfetsch~\cite{fischerpfetsch2015}
emphasize cuts and branching techniques for problems with overlapping SOS1 constraints,
that is, sets of complementarity constraints that have variables in common;
this structure can be represented with conflict graphs and can be exploited
in the derivation of valid cutting planes and in the construction of sophisticated branching rules
building on the ideas of Beale and Tomlin~\cite{bealetomlin1970}.
In our formulation, each variable appears in at most one complementarity
constraint, so the nice techniques of Fischer and Pfetsch would not be helpful.
 
\subsection{LPCC using Branch-and-Cut\label{sec:LPCC using branch-and-cut}}

In this paper, we propose a branch-and-cut algorithm for solving the general LPCC problem (\ref{eq:general LPCC}).
In \S\ref{sec: Preprocessing Phase}, we  describe the preprocessing phase of our algorithm:
in \S\ref{sec: LPCC Feasibility Recovery}, a heuristic feasibility recovery procedure is developed to recover a feasible solution of the LPCC which provides a valid upper bound of the LPCC;
and in \S\ref{sec: Cutting Planes Generation and Selection}, the strategy of generating and selecting from various types of cutting planes we studied in~\cite{JEMPYuMOPTA} is discussed, which could sharpen the LP relaxation and improve the initial lower bound of LPCC.
In \S\ref{sec: Branch-and-Bound Phase}, we present the second phase of our algorithm: branch-and-bound.
Various node selection strategies and branching complementarity selection strategies are discussed in \S\ref{sec: Branching Complementarity Selection} and \S\ref{sec: Node Selection}.
Our proposed algorithm is able to characterize infeasible and unbounded LPCC problems as well as solve problems with finite optimal value.
The algorithm is summarized in~\S\ref{sec: General Scheme of the Algorithm}.
In \S\ref{sec: Computational Result for branch-and-cut algorithm}, we show the computational results of our branch-and-cut algorithm on solving randomly generated LPCC instances.

In the  MIP formulation (\ref{eq:MIP LPCC}), the binary vector $z$ is only used to model the complementary relationship of the LPCC, and except for the complementarity constraints it does not interact with $x$ and $y$ at all. This observation motivates us to enforce the complementarities through 
a specialized branching scheme, i.e., branch on complementarities directly without introducing the binary vector $z$. This kind of specialized branching approach has been studied to solve several problems such as generalized assignment problems~\cite{farias1}, nonconvex quadratic programs~\cite{vandenbussche2}, nonconvex piecewise linear optimization problems~\cite{keha2},
and problems with overlapping SOS1 constraints~\cite{farias8,fischerpfetsch2015}.
The obvious advantage of using a specialized branching approach for solving the LPCC is that we no longer need $\theta$ in the formulation, and therefore this approach is also applicable for the case when $y$ or $w$ is unbounded.
In fact, even if we know such a $\theta$ exists, the cost of computing a valid $\theta$ could be very expensive especially when $m$ is very large.
Moreover, introducing the binary vector $z$ will lead to an increase in both the number of variables and the number of constraints,
and these Big-M type constraints are usually not tight which will lead to a number of violated complementarities in the solution of the relaxation.

\section{Preprocessing Phase \label{sec: Preprocessing Phase} }

When the initial LP relaxation is bounded below, the preprocessing phase will be invoked,
consisting of a feasibility recovery process and a cutting plane selection and management process.
The feasibility recovery process may provide a valid upper bound for the LPCC, while the cutting plane selection and management process may provide a better lower bound for the LPCC. Both processes may provide a good starting point for the second phase of our algorithm: branch-and-bound. 

\subsection{LPCC Feasibility Recovery \label{sec: LPCC Feasibility Recovery} }
Finding a good feasible solution to an LPCC is an essential component of our branch-and-cut algorithm for globally resolving the LPCC. A good upper bound can help prune nodes quickly, and avoid unnecessary branching. Notice that here we assume the initial LP relaxation is bounded when we apply our feasibility recovery procedures.
Our feasibility recovery procedures have some similarities
to feasibility pumps for MIP and MINLP~\cite{fischetti11}.

For ease of discussion, we first introduce some notation and definitions. 
\begin{definition}
Given any binary vector $z$ with dimension $m$, we define the linear program LPCC($z$) as follows:
\begin{equation} \label{eq: a piece of LPCC}
\begin{array}{lll}
\displaystyle{
{\operatornamewithlimits{\mbox{minimize}}_{(x,y)}}
} & c^Tx + d^Ty\\ 
\mbox{subject to} & Ax + By  \geq  b \\ 
& x \, \geq \, 0  \\
 & 0  \leq  y \\
 & 0 \leq  q + Nx + My \\
 & 0 \geq y_i & \mbox{ if } z_i=0 \\
 & 0  \geq (q + Nx + My)_i & \mbox{ if }  z_i=1.
\end{array} \end{equation}
LPCC($z$) is a so-called piece of the LPCC corresponding to the binary vector $z$. 
\end{definition} 
\begin{definition}
The feasibility gap of the piece of an LPCC corresponding to the binary vector~$z$, denoted by FG-LPCC($z$), is the optimal value of the following linear program:
\begin{equation} \label{eq: feasible level of a piece of LPCC}
\begin{array}{ll}
\displaystyle{
{\operatornamewithlimits{\mbox{minimize}}_{(x,y,w)}}
} & (\textbf{1}-z)^Ty + z^Tw\\ 
\mbox{subject to} & Ax + By  \geq  b \\ 
& x \, \geq \, 0  \\
 & 0  \leq  y \\
 & 0 \leq w:= q + Nx + My  \\
\end{array} \end{equation}
where \textbf{1} is the vector with all  components equal to 1.
\end{definition}
Based on the above two definitions, it is obvious that the following proposition is true:
\begin{proposition}
$LPCC(z)$ is feasible \textit{if} and \textit{only if} FG-LPCC($z$)=0
\end{proposition}
\begin{definition}
Binary vectors $z$ and $z'$ are {\em adjacent}
if there is exactly one component that is different between $z$ and $z'$.
\end{definition}
\begin{definition}
If binary vectors $z$ and $z'$ are adjacent and FG-LPCC($z$) $<$ FG-LPCC($z'$), then $\Delta z=z-z'$ is a
{\em feasibility gap descent direction} for $z'$.  
\end{definition}

Just like mixed integer programs, it is often
a good idea to recover a feasible solution based on the LP relaxation solution. The most intuitive recovery process is to round the LP relaxation solution into a solution that satisfies all the complementarity constraints.
We use this rounding procedure to initialize a new local search feasibility recovery process,
detailed in 
Procedure~\ref{proc: Local search feasibility recovery process}.
Notice that we define search $breadth$ as the number of candidates that we are going to select from binary vectors which are adjacent to the initial~$z^*$,
and search $depth$ as the maximum number of iterations that we are going to perform for each candidate. We can set search $breadth$ and search $depth$ to control the local search process.
The proposed local search procedure can be used to find a feasible solution, 
although the quality of the recovered feasible solution is not guaranteed.
We use optimality based bound tightening \cite{gleixner16,maranas1997,zamora1999}
to resolve this issue,
refining the local search feasibility recovery procedure through the
addition of the constraints $lb_{search}\le c^Tx+d^Ty\le ub_{search}$
to~(\ref{eq: feasible level of a piece of LPCC}) when computing the feasibility gap.
Procedure~\ref{proc: Refined local search feasibility recovery process} describes this refined feasibility recovery procedure. We will demonstrate the computational results of our proposed local search feasibility recovery process in \S\ref{sec: Computational Result for branch-and-cut algorithm}.
See Fischer and Pfetsch~\cite{fischerpfetsch2015} for primal
heuristics that can be used when a variable appears in more than
one complementarity constraint.

\begin{algorithm}
\SetAlgorithmName{Procedure}{}

\SetKwInOut{Input}{input}\SetKwInOut{Output}{output}

\Input{the LP relaxation solution of the original LPCC: $x^*$, $y^*$, $w^*$, search depth parameter $depth$, search breadth parameter $breadth$}
\Output{recovered feasible LPCC solution or failed to recover the solution}
\BlankLine
Initialization:  Set binary vector $z^*=0$\;

\For{$i \leftarrow 1$ \KwTo $m$}
{
	\lIf{$y^*_i<w^*_i$}{$z^*_i$=0\;}
	\lElse{$z^*_i$=1\;}
}
Solve (\ref{eq: feasible level of a piece of LPCC}) to get \textit{ FG-LPCC($z^*$)}\;
\If{FG-LPCC($z^*$)==0}{solve \textit{LPCC($z^*$)}, and \textbf{return} the optimal solution to \textit{LPCC($z^*$)}\;}
\Else{ let $A(z^*)$ denote the set of binary vectors that are adjacent to $z^*$\;
	\ForEach{$z \in A(z^*)$}{solve (\ref{eq: feasible level of a piece of LPCC}) to get \textit{FG-LPCC($z$)}\;
	insert $z$ into a sorted queue $Q$ with nondecreasing order on \textit{FG-LPCC($z$)}\;}
	Let $r_b=0$\;
	\While{$Q$ is not empty and $r_b\le breadth$}
	{	
		$r_b=r_b+1$\;
		pop the top element $\bar{z}$ in $Q$, and delete this element from $Q$\;
		let $z=\bar{z}$ and $r_d=0$\;
		\While{there exists any feasibility gap descent direction $\Delta z$ for $z$ and $r_d\le depth$ }
		{	
			pick a feasibility gap descent direction $\Delta z$\;
			$z=z+\Delta z$\;
			$r_d=r_d+1$\;
		}
		\If{\textit{FG-LPCC($z$})==0}{solve \textit{LPCC($z$)}, and \textbf{return} the optimal solution to \textit{LPCC($z$)}\;}
	}
}
\textbf{return} feasibility recovery failed\;
	
\caption{Local search feasibility recovery process}\label{proc: Local search feasibility recovery process}
\end{algorithm} 

\begin{algorithm}
\SetAlgorithmName{Procedure}{}

\SetKwInOut{Input}{input}\SetKwInOut{Output}{output}

\Input{the known valid upper bound of LPCC $ub_{initial}$, parameter  $searchGap_{min}$}
\Output{refined feasible LPCC solution or failed to refine the known feasible solution}
\BlankLine
Initialization:  Set $lb_{search}=$optimal value of the LP relaxation of LPCC and $ub_{search}=ub_{initial}$;
				add $lb_{search}\le c^Tx+d^Ty\le ub_{search}$ into (\ref{eq: feasible level of a piece of LPCC})\;
\While{$ub_{search}-lb_{search}>searchGap_{min}$}
{
	solve LP relaxation of LPCC with constraints $lb_{search}\le c^Tx+d^Ty\le ub_{search}$ \;
	apply Procedure \ref{proc: Local search feasibility recovery process} to recover a feasible solution\;
	\If{recovery process succeed}{
	update the refined feasible solution with recovered solution\;	
	update $ub_{search}$ with the newly recovered solution\; 
	$ub_{search}=(lb_{search}+ub_{search})/2$\;}
	\Else{$lb_{search}=(lb_{search}+ub_{search})/2$\;}
}
\If{refined feasible solution has been updated}{\textbf{return} refined feasible solution}
\Else {\textbf{return} feasibility refinement failed} 
\caption{Refined local search feasibility recovery process}\label{proc: Refined local search feasibility recovery process}
\end{algorithm}

\subsection{Cutting Plane Generation and Selection\label{sec: Cutting Planes Generation and Selection} }
The second key step in our preprocessing phase is the generation and selection of cutting planes.
We have discussed various valid linear constraints and second order cone constraints that can be used to tighten the initial relaxation of LPCC in \cite{JEMPYuMOPTA},
and have shown the computational results of these valid constraints individually. 
As important as finding these cutting planes is the selection of the cuts that actually should be
added to the initial LP relaxation. In this section, we will describe our detailed procedure to generate and select our cutting planes.
Note that we will only add cutting planes at the root node, and perform the generation of each type of cut in rounds and in the following order:
\begin{itemize}
\item
Disjunctive cuts and Simple cuts
\item
Bound cuts
\item
Linear cuts derived from second order cone constraints
\end{itemize}

We use the computational results with these cutting planes in \cite{JEMPYuMOPTA} to guide the cut generation process. The details of generation and selection rule are described as follows.

\subsubsection{Disjunctive cuts and simple cuts}
These cuts exploit the disjunctive constraints: for each~$i$, either $y_i \leq 0$ or $w_i \leq 0$.
The solution to the LP relaxation typically violates a number of these disjunctions,
and disjunctive cuts can either be generated by solving a supplemental linear program,
or by examining the optimal tableau for the LP relaxation.
Based on our computational experience, it seems that general disjunctive cuts and simple cuts~\cite{balas23,audet5,audet6}
are the weakest cuts among our three type of cutting planes, but they are the cheapest to generate.
Therefore we generate this type of cut first.
The solving time of CPLEX for our test instances became  worse when we added all of the generated disjunctive cuts or simple cuts to the root node even though the value of the initial LP relaxation was improved by these cuts,
because the initial LP became too large.
Moreover, due also to the high cost of generating general disjunctive cuts, we only generate $\lfloor m/100 \rfloor$
rounds of general disjunctive cuts and for each round we only generate at most 3 general disjunctive cuts instead of generating disjunctive cuts for each violated complementarity constraint.

The values of $y_iw_i$ in the optimal solution to the LP relaxation
are sorted in nonascending order and we select complementarity constraints
with index that corresponds to the largest three products.
After each round of generating cuts, we will remove every cut whose corresponding slack variable is basic in the relaxed LP, in order is to keep the size of the relaxed LP small.
After generating the general disjunctive cuts, $\lfloor m/10 \rfloor$ rounds of simple cuts will be added.
Since a simple cut is derived from the simplex tableau with almost no cost,
we will generate simple cuts for every violated complementarity constraint in each round,
and also remove every cut whose corresponding slack variable is basic in the relaxed LP
after each round of generating cuts.

\subsubsection{Bound cuts}
Upper bounds $u^y_i$ and $u^w_i$ on $y_i$ and $w_i$ can be used in the bound cut
\begin{equation}
u^w_i y_i + u^y_i w_i \, \leq \, u^w_i u^y_i
\end{equation}
for any pair of complementary variables $y_i$ and~$w_i$.
Strengthening the upper bounds seems
very important for the branch-and-bound routine of CPLEX for solving our instances,
and the bound cuts also improve the initial lower bound dramatically.
However, the major drawback of bound cuts is that they are very expensive to generate, especially when $m$, the number of complementarity constraints, is very large. Therefore, we will only compute bounds for at most 5 pairs of complementary variables, and the selection of these complementary variables is the same as the selection of complementarity constraints to generate disjunctive cuts.
An upper bound $u^y_i$ for $y_i$ can be found by solving the linear program
\begin{equation} \label{eq:uyi}
u^y_i  \, = \, \begin{array}[t]{ll}
\displaystyle{
{\operatornamewithlimits{\mbox{maximize}}_{(x,y,w)}}
} & y_i   \\
\mbox{subject to} & Ax + By  \geq  b \\ 
& x \, \geq \, 0  \\
 & 0  \leq  y \leq u^y \\
 & 0 \leq  q + Nx + My \, = \, w \,  \leq  u^w  \\
 & c^Tx + d^Ty \leq ub  \\
 & u^w_j y_j + u^y_j w_j \, \leq \, u^w_j u^y_j  \quad  \forall \, j \mbox{ with known bounds } u^y_j, u^w_j
\end{array} \end{equation}
where $ub$ is a known upper bound on the optimal value of the LPCC.
A similar LP can be constructed to get bounds on~$w$.

We also investigated improving the bound cuts by splitting the variables.
In particular, two versions of problem (\ref{eq:uyi}) could be solved,
one with the additional constraint $y_k=0$ and the other with the additional constraint $w_k=0$,
for some index~$k \neq i$.
The maximum of the optimal values of these two problems could potentially improve on the initial upper bound.
For our test instances, the additional computational work involved in computing these improved
bounds did not improve the overall computational time, so this splitting is not included in our results.

\subsubsection{Linear cuts from second order cone constraints}
Based on the computational results of \cite{JEMPYuMOPTA},  cuts derived from a certain second order cone constraint can significantly improve the initial lower bound of our  instances with relatively low generating cost compared to bound cuts when $n << m$, provided $M$ is positive semidefinite.
These cuts arise from linearizing the term $y^TNx$, using
McCormick inequalities~\cite{GPMccormick2} to tighten the linearization,
and handling the $y^TMy$ term appropriately.
Details can be found in~\cite{JEMPYu_QCQP}.
The constraints can be tightened by refining bounds.
We did not use these cuts in the computational results reported in this paper, because
of difficulties with ensuring $M$ was regarded as numerically positive semidefinite
by CPLEX.

\subsection{Overall Flow of the Preprocessor\label{sec:Overall Flow of the Preprocessor}}

The preprocessor consists of the following steps:
\begin{enumerate}
\item Apply the feasibility recovery routine to recover a feasible solution.
\item Generate $\lfloor m/100 \rfloor$ rounds of general disjunctive cuts.
\item Generate $\lfloor m/10 \rfloor$ rounds of simple cuts.
\item Apply 4 bound refinements
and generate bound cuts.
\end{enumerate}
\setcounter{remark}{0}

We apply Procedures \ref{proc: Local search feasibility recovery process}
and~\ref{proc: Refined local search feasibility recovery process}
as the default feasibility recovery procedure due to run time considerations. Other feasibility recovery procedures and refinements can also be invoked if required for solving special classes of problems. 
The number of rounds for generating each type of cutting plane can be modified by changing the parameter settings.
The current setting is based on the computational experience in~\cite{JEMPYuMOPTA}.

An additional preprocessing procedure undertaken at each node is
the complementary variable fixing process,
which is detailed in~\S\ref{sec: Node Pre-solving and warm start}.

\section{Branch-and-Bound Phase \label{sec: Branch-and-Bound Phase} }
The branch-and-bound routine needs to be invoked to solve the problem exactly if the initial LP relaxation is unbounded or the preprocessing phase is unable to close $100\%$ of the gap for the bounded case.
The branching is imposed on the complementarity constraint directly, and two subproblems (nodes) will be generated by enforcing either side of the
pair of complementary variables to its lower bound zero. Just like a
branch-and-bound based MIP solver, there are two key ingredients in our branch-and-bound routine:  branching complementarity selection and node selection.
Branching complementarity selection is the procedure to select the complementarity constraint to be  branched on, and it is the same as the ``variable selection" in mixed integer programming.
In \S\ref{sec: Branching Complementarity Selection},
we  present our branching strategy which is based on the ideas of  three classic branching rules and also some new proposed ideas designed for the LPCC problem. Node selection is the procedure to select the next subproblem from the node tree to be processed. In \S\ref{sec: Node Selection}, we will present and compare different node selection strategies. Besides these two key ingredients, in \S\ref{sec: Node Pre-solving and warm start} we will describe the node pre-solving procedure used in our algorithm to pre-process the nodes during the branch-and-bound process. The general branch-and-bound routine for handling the bounded case and unbounded case of LPCC are described in \S\ref{sec: Overall Flow of Branch-and-Bound for LPCC Bounded Case} and
\S\ref{sec: Overall Flow of Branch-and-Bound for LPCC Unbounded Case} respectively.

\subsection{Branching Complementarity Selection \label{sec: Branching Complementarity Selection} }
The branching rule is the key ingredient of any branch-and-bound algorithm. Good branching strategies are extremely important in practice for solving mixed integer programs, although currently there is no existing theoretical best branching strategy.
We will first present three classic branching strategies for solving mixed integer programs that have been studied in the literature. The reader can  refer to Linderoth and Savelsbergh~\cite{linderoth1},
F\"ugenschuh and Martin~\cite{fugenschuh1} and Achterberg et al~\cite{achterberg1} for a comprehensive study of branch-and-bound strategies for mixed integer programming.
We will present our branching strategy based on the ideas of these branching strategies. The computational results that compare various branching strategies will be shown in \S\ref{sec: Computational Result for branch-and-cut algorithm}.

We first give some definitions related to our branching routine for the LPCC problem. For easy discussion, if the LP relaxation of the LPCC is unbounded below, we represent its lower bound as $-\infty$.
Suppose that we have an LPCC problem $Q$ and the set $I$ is the index set of complementarity constraints. If the current solution to the LP relaxation of $Q$ is not a feasible solution to LPCC
(for the unbounded case,we consider an unbounded ray of the LP relaxation instead of solution to
the LP relaxation),
then we can pick an index $i\in I$ with $y_iw_i>0$ and obtain two subproblems (nodes): one by adding the constraint $y_i\leq0$ (named the left child node, denoted by $Q^y_i$) and one by adding the constraint $w_i\leq0$ (named the right child node, denoted by $Q^w_i$).
We refer to this as branching on complementarity~$i$.
For the bounded case, if we denote the objective value of the LP relaxation of $Q$ as $c_Q$ and the objective value of the LP relaxation of its two child nodes as $c_{Q_i^y}$ and $c_{Q_i^w}$ respectively, then the objective value changes caused by branching on the $i$th
complementarity are $\Delta_i^y=c_{Q_i^y}-c_Q$ and $\Delta_i^w=c_{Q_i^w}-c_Q$. We usually use the improvement of objective value of the LP relaxation to measure the quality of  branching on the $i$th complementarity.
Our implementation supports fixing multiple complementarity constraints at one time, but by default we will only select to branch on one complementarity.
Based on the results of testing our  instances and the computational results of solving various MIP problems in the literature, multiple way branching is rarely better than two way branching. 

The generic procedure for selecting the branching complementarity can be described in
Procedure~\ref{proc: generic complementarity selection procedure}.
\begin{algorithm}
\SetAlgorithmName{Procedure}{}

\SetKwInOut{Input}{input}\SetKwInOut{Output}{output}

\Input{the LP relaxation solution of the current processing node $Q$ or the unbounded ray to the LP relaxation if the LP relaxation is unbounded: $x^*$, $y^*$, $w^*$}
\Output{the selected branching index $i\in I$ of a complementarity constraint}
\BlankLine
\begin{enumerate}
\item Let $\tilde{I}=\{j\in I \, | \, y^*_j \,  w^*_j \, > \, 0\}$ denote the index set of violated \\ complementarity constraints.
\item Compute a branching score $s_j \in \R^+ $ for all candidates $j\in \tilde{I}$.
\item Select the selected branching index $i\in \tilde{I}$ with $s_i=\max_{k\in \tilde{I}}\{s_k\}$.
\end{enumerate}
Return selected branching index $i$.
\caption{Generic complementarity selection procedure}\label{proc: generic complementarity selection procedure}
\end{algorithm} 
The score function in Step 2 of this procedure needs to evaluate the two child nodes that could be generated by the branching, and map these two effectiveness values onto a single score value.
Different choices for the effectiveness values are given later.
Suppose $q^y$ and $q^w$ are the effectiveness values of the two child nodes generated by a branching. In the literature, the score function usually has one of the following forms:
\begin{equation} \label{eq: weighted sum branch score function}
score(q^y,q^w)=(1-\mu)\cdot \min\{q^y,q^w\}+\mu\cdot \max\{q^y,q^w\}
\end{equation}
or
\begin{equation} \label{eq: product form branch score function}
score(q^y,q^w)=\max\{q^y,\epsilon\}\cdot \max\{q^w,\epsilon\}
\end{equation}
Here $\mu$ is a number between 0 and 1, and it is usually an empirically determined constant or a dynamic parameter adjusted through the course of branching process.
We chose $\epsilon=10^{-6}$ to enable the comparison when either $q^y$ or $q^w$ is zero. Based on the computational experience in \cite{Achterberg.thesis}, the product form is superior to the weighted sum form for solving MIP problems. Therefore, in our algorithm, we chose to use the product form to map the effectiveness values from two child nodes onto a single value. 

In the following we will present three classic branching strategies for solving an MIP in terms of our branching on complementarity scheme: \textit{Strong Branching} (apparently originally developed in the work leading
up to~\cite{applegate6}), \textit{Pseudocost Branching}~\cite{benichou1970}
and \textit{Inference Branching}~\cite{Achterberg.thesis}.
In fact, all of these branching routines are just variants of  Procedure \ref{proc: generic complementarity selection procedure} with  different score functions.

\subsubsection{Strong branching}
The idea of \textit{Strong Branching}~\cite{applegate6}
is to test the branching candidates by temporarily enforcing either side of a complementarity constraint and solving the resulting LP relaxation to a certain level, then select the one that can lead to the largest lower bound improvement.
\textit{Full Strong Branching} will compute $\Delta_i^y$ and $\Delta_i^w$
for each branching complementarity candidate $i\in \tilde{I}$, and use the score($\Delta_i^y$, $\Delta_i^w$) as the effectiveness values
in the form of either (\ref{eq: weighted sum branch score function}) or (\ref{eq: product form branch score function}) as its score function. \textit{Full Strong Branching} can be seen as the locally best branching strategy in terms of lower bound improvement. However the computational cost of \textit{Full Strong Branching} is very high, since in order to evaluate the score function for each complementarity candidate, we need to solve two resulting LP relaxations to optimality. There are usually two ways to speed up \textit{Full Strong Branching}: one is to only test a subset of the candidate set instead of considering all the candidates, and another is to perform a limited number of simplex iterations and estimate the objective value change based on that. In our branch-and-bound algorithm, we have implemented the \textit{Full Strong Branching} routine, and also we adopt the former idea to speed up the \textit{Full Strong Branching}: as long as the objective value of LP relaxation of either side of the child nodes hits some threshold, we will select this branching candidate and exit the selection routine; we set the median value of the lower bound of unsolved nodes in the current search tree as this threshold.

A version of strong branching was used by Fischer and Pfetsch~\cite{fischerpfetsch2015}
in their branch-and-cut approach for problems with overlapping SOS1 constraints.

\subsubsection{Pseudocost branching   \label{subsec:pseudocost}}
\textit{Pseudocost Branching}~\cite{benichou1970}
uses the branching history to estimate the two objective changes of the child nodes without actually solving them. In other words, \textit{Pseudocost Branching} is a branching rule based on the historical performance of complementarity branching on complementarities which have already been branched. Let $\varsigma_i^y$ and $\varsigma_i^w$ be the objective gain per unit change at node $Q$ after branching on complementarity $i$ by enforcing $y_i$ or $w_i$ to zero, that is
\begin{equation}
\varsigma_i^y=\dfrac{\Delta_i^y}{y_i^*} \mbox{  and  } \varsigma_i^w=\dfrac{\Delta_i^w}{w_i^*}
\end{equation}
where $y_i^*$ and $w_i^*$ are the violation of complementarity $i$ corresponding to the LP relaxation solution of~$Q$. Let $\sigma_i^y$ denote the sum of  $\varsigma_i^y$ over all the processed nodes
where complementarity $i$ has been selected as the branching complementarity and resulting child node $Q_i^y$ has been solved and was feasible. Let $\eta_i^y$ denote the number of these problems, and define $\sigma_i^w$and $\eta_i^w$ in the same way for the other side of the complementarity. Then the pseudocost of branching on complementarity $i$ can be calculated as the arithmetic mean of objective gain per unit change:
\begin{equation}
\Psi_i^y=\dfrac{\sigma_i^y}{\eta_i^y} \mbox{  and  } \Psi_i^w=\dfrac{\sigma_i^w}{\eta_i^w}
\end{equation}
Therefore given the violated complementarity $i$ corresponding to the LP relaxation of $Q$, it is reasonable to use $\Psi_i^y\cdot y_i^*$ and $\Psi_i^w\cdot w_i^*$ to estimate $\Delta_i^y$ and $\Delta_i^w$ respectively. We call the branching rule that uses
the score function $score(\Psi_i^y\cdot y_i^*,\Psi_i^w\cdot w_i^*)$ in step 2 of Procedure \ref{proc: generic complementarity selection procedure} as \textit{Pseudocost Branching}. Notice that at the beginning of the branch-and-bound procedure, the pseudocost is uninitialized for all the complementarities. One way to handle a complementarity with an uninitialized pseudocost is to replace its pseudocost with the average of the pseudocosts of the complementaries whose pseudocosts have been initialized, and set the pseudocost as 1 if all the complementarities are uninitialized. Applying strong branching to the nodes whose tree depth level is less than a given level is another way to initialize the pseudocosts.
More recently, Achterberg et al~\cite{achterberg1} proposed a more general pseudocost initialization method, and named the corresponding branching rule as \textit{Reliability Branching}. In our implementation, we include the pseudocost as part of our branching score, and we choose to apply strong branching to nodes whose tree depth level is less than 7 to initialize the pseudocost.

\subsubsection{Inference branching}
The branching decision of strong branching and pseudocost branching are both based on the change of objective value of the LP relaxation, while \textit{Inference Branching}~\cite{Achterberg.thesis}
is quite different from the above two branching strategies.
\textit{Inference Branching} checks the impact of branching on changing the bounds of other variables.
As with pseudocosts,
historical information is typically used to estimate the deductions on bounds of the variables,
and the inference value can be calculated as the arithmetic mean of the number of bound deductions.
The domain propagation process is a node pre-solving process to detect the bound change of the variables
and is discussed in~\S\ref{sec: Node Pre-solving and warm start}.
In our implementation, we use a similar idea to inference branching: instead of evaluating
the inference value, we estimate the \textit{complementarity satisfaction level} after branching on a complementarity,
leading to the quantity~$s^{SL}_i$ below.

\subsubsection{Hybrid branching strategy for the LPCC (bounded case)  \label{subsec:branch_bdd}}

Our branching strategy for the bounded case combines the ideas of the above three classic branching strategies, and additionally we also include some new score values into our branching score function which are specialized for the LPCC problem.

In our implementation, the default branching strategy will apply the full strong branching strategy for the nodes whose depth level are no larger then $7$. The reason for doing that is because it is usually quite important to make the right branching decision at the beginning, and also we can use strong branching to initialize the pseudocosts and another score value that we will propose next. For the nodes whose tree depth are larger than $7$, we will use a weighted sum formula to combine four score values for each violated complementarity. Among these four score values, two of them are only based on the current node~$Q$, and the other two are based on historical branching information.
For the violated complementarity~$i$, these four score values are listed as follows:
\begin{enumerate}
\item $s_i^{VL}$: score of Violation Level. Suppose $y_i^*$ and $w_i^*$ are the violation of complementarity~$i$ corresponding to the LP relaxation of~$Q$, then we define
$$
s_i^{VL}=\sqrt{y_i^*\cdot w_i^*}
$$
\item $s_i^{ED}$: score of Euclidean Distances from the LP relaxation solution of Q to the two hyperplanes corresponding to $y_i=0$ and $w_i=0$. Recall that since $y_i^*\cdot w_i^*>0$, we can represent the complementary variables $y_i$ and $w_i$ with the non-basic variables in the optimal simplex tableau of $Q$
\begin{eqnarray}
y_i&=&y_i^*-\displaystyle{\operatornamewithlimits{\sum}_{j\in NB}}a^{yi}_j\xi_j
\label{eq:y nonbasic}\\
w_i&=&w_i^*-\displaystyle{\operatornamewithlimits{\sum}_{j\in NB}}a^{wi}_j\xi_j
\label{eq:w nonbasic}
\end{eqnarray}

We use the Euclidean distance from the LP relaxation solution to the two hyperplanes
$$
\displaystyle{\operatornamewithlimits{\sum}_{j\in NB}}a^{yi}_j\xi_j = y_i^* \mbox{ and }
\displaystyle{\operatornamewithlimits{\sum}_{j\in NB}}a^{wi}_j\xi_j = w_i^*
$$
to define $s_i^{ED}$ as follows:
$$
s_i^{ED}=\sqrt{\dfrac{y_i^*\cdot w_i^*}{\sqrt{\Vert a^{yi}\Vert\cdot\Vert a^{wi}\Vert}}}
$$
\item $s^{PC}_i$: score of Pseudo Cost. We use
the following small modification to the pseudcost calculation of~\S\ref{subsec:pseudocost}:
$$
s^{PC}_i=\sqrt{\max\{\Psi_i^y\cdot y_i^*,\epsilon\}\cdot \max\{\Psi_i^w\cdot w_i^*,\epsilon\}}
$$
\item $s^{SL}_i$: score of complementarity Satisfaction Level. We define the complementarity satisfaction level as the proportion of the satisfied complementarities corresponding to the LP relaxation solution of the child node after branching. Intuitively we want to select a branching complementarity that will lead to more satisfied complementarities. To estimate this complementarity satisfaction level, we collected the historical information to compute the average complementarity satisfaction level for both sides of the complementarity 
$$
\Phi_i^y=\dfrac{\varphi_i^y}{\eta_i^y} \mbox{ and } \Phi_i^w=\dfrac{\varphi_i^w}{\eta_i^w}
$$
Here $\varphi_i^y$ is the sum of the proportion of complementarity satisfaction levels over all the prior nodes, where complementarity $i$ has been selected as the branching complementarity, and $\eta_i^y$ is the total number of these nodes. We define $\varphi_i^w$ and $\eta_i^w$ to be the analogous value for the other side of complementarity. Then the score of the complementarity Satisfaction Level can be calculated as
$$
s^{SL}_i=\sqrt{\Phi_i^y\cdot\Phi_i^w}
$$
\end{enumerate}
We scale the score vectors using their 2-norms, and the following formula is the
branching score function that we used to evaluate the score for each violated complementarity:
\begin{equation} \label{eq: score function of branching complementarity selection}
s_i=\omega^{VL}\left( \dfrac{s_i^{VL}}{\Vert s^{VL} \Vert}\right) +\omega^{ED}\left( \dfrac{s_i^{ED}}{\Vert s^{ED} \Vert}\right) +\omega^{PC}\left( \dfrac{s_i^{PC}}{\Vert s^{PC} \Vert}\right) +\omega^{SL}\left( \dfrac{s_i^{SL}}{\Vert s^{SL} \Vert}\right) 
\end{equation}
By default, the weight is set as $\omega^{VL}=1$, $\omega^{ED}=0.5$, $\omega^{PC}=0.25$ and $\omega^{SL}=0.5$. Note that  setting different weights for each score value will lead to different branching behaviour. In 
\S\ref{sec: Computational Result of branch-and-cut}, we will show the computational results of solving our  LPCC instances with different weights of the score value.

\subsubsection{Hybrid branching strategy for LPCC (unbounded case)  \label{subsec:branch_unbdd}}

Our branching strategy for the unbounded case is slightly simpler than the one for the bounded case.
We will still apply full strong branching to the nodes whose tree depth level is no larger than $7$.
However, for the remaining unbounded nodes we will only use $s_i^{VL}$ as the branching score to make the branching decision.

\subsection{Node Selection \label{sec: Node Selection} }
In addition to selecting which complementarity to branch on, another question is which subproblem (node) we should pick to process. There are two major criteria for selecting the next subproblem to be processed. 
\begin{enumerate}
\item finding feasible LPCC solutions to improve the upper bound of the LPCC problem which leads to pruning the nodes by bounding, leading to a Depth First Search strategy. 
\item improving the lower bound as fast as possible,
leading to a Best-Bound strategy.
\end{enumerate}

In our implementation of the branch-and-bound routine we use a \textbf{\textit{Best-Bound}} strategy to select the next node to be processed, since we want to solve the problem to optimality as fast as possible. Notice that for the \textit{Best-Bound}, it is possible that there are several nodes with the same lower bound. For that case, we will select the most recently generated node as the next node to be processed.

\subsection{Node Pre-solving
\label{sec: Node Pre-solving and warm start} }

The major task of our node pre-solving procedure is to tighten the domains of complementary variables $y_i$ and $w_i$ and try to fix the complementary variables.
In order to facilitate the discussion, here we can assume that each $\xi_i$ in
(\ref{eq:y nonbasic}) and (\ref{eq:w nonbasic})
is a non-negative variable with zero lower bound. Therefore we have the following result: if $a^{yi}_j\leq 0, \forall j\in NB$, then we have $y_i\geq \hat{y_i}$, and therefore $w_i=0$; if $a^{wi}_j\leq 0, \forall j\in NB$, then we have $w_i\geq \hat{w_i}$, and therefore $y_i=0$.
This complementary variable fixing check is performed before we branch on the complementarity constraint.

\section{General Scheme of the Branch-and-Cut Algorithm for Solving LPCC\label{sec: General Scheme of the Algorithm} }

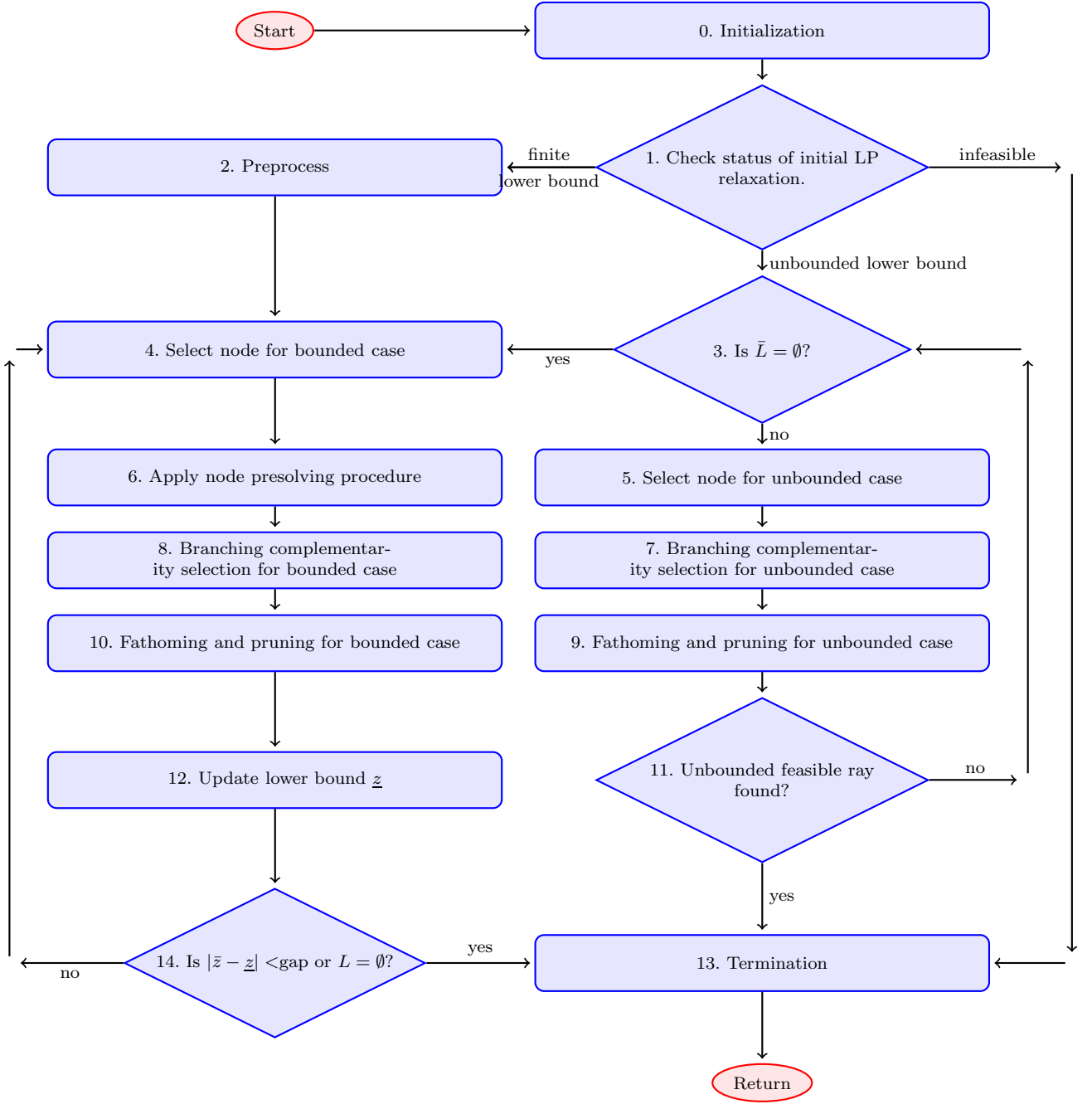
\begin{figure}
\begin{center}
\begin{tikzpicture}
  [auto,
decision/.style={shape aspect=2, diamond, draw=blue, thick, fill=blue!10, text width=4cm,align=flush center,
inner sep=1pt,  font=\footnotesize},
          block/.style
={rectangle, draw=blue, thick, fill=blue!10,
text width=7cm,align=center, rounded corners,
minimum height=3em,  font=\footnotesize},
line/.style
={draw, thick, ->, 
shorten >=2pt},
cloud/.style
={draw=red, thick, ellipse,fill=red!10,
minimum height=2em,  font=\footnotesize}]

                \matrix [column sep=5mm,row sep=4mm] {

&     \node [cloud] (start) {Start};   & \node [block] (init) {0. Initialization};  &
 \\
  &   \node [block] (leftpt) {2. Preprocess}; 
  &  \node [decision] (status) {1. Check  status of  initial LP relaxation.};
  & & \node (rightpt) {};  \\ 
   \node (leftpt2) {}; & \node [block] (bddnode) {4. Select node for bounded case};
 &  \node [decision] (Lbarempty) {3. Is $\bar{L}=\emptyset$?};  & \node (rightpt2) {};  \\
   & \node [block] (presolve) {6. Apply node presolving procedure};
  & \node [block] (unbddnode) {5. Select node for unbounded case};  \\
   &\node [block] (bddbranch) {8. Branching complementarity selection for bounded case};
   &\node [block] (unbddbranch) {7. Branching complementarity selection for unbounded case};  \\
   &  \node [block] (bddfathom) {10. Fathoming and pruning for bounded case};
   &  \node [block] (unbddfathom) {9. Fathoming and pruning for unbounded case}; \\
   &  \node [block] (updatelb2) {12. Update lower bound $\underline{z}$};
   & \node [decision] (ray) {11. Unbounded feasible ray found?};
   & \node (rightpt3) {};  \\
   \node (leftpt3) {};  &
   \node [decision] (stopbdd) {14. Is $|\bar{z}-\underline{z}|<$gap or $L=\emptyset$?};
    & \node [block] (term) {13. Termination};
    && \node (rightpt4) {}; \\
    && \node [cloud] (return) {Return};  \\
};
\begin{scope}[every path/.style=line]
    \path  (start) -- (init);
    \path  (init) -- (status);
    \path  (status) -- node [midway, above] {\footnotesize finite} (leftpt);
    \path  (status) -- node [midway, below] {\footnotesize lower bound} (leftpt);
    \path (status) -- node [midway] {\footnotesize infeasible} (rightpt);
    \path (rightpt) -- (rightpt4);
    \path (status) -- node [midway] {\footnotesize unbounded lower bound} (Lbarempty);
    \path (leftpt) --  (bddnode);
    \path  (Lbarempty) -- node [midway] {\footnotesize yes} (bddnode);
    \path (Lbarempty) -- node [midway] {\footnotesize no} (unbddnode);
    \path (leftpt2) -- (bddnode);
    \path (bddnode) -- (presolve);
    \path  (unbddnode) -- (unbddbranch);
    \path (presolve) -- (bddbranch);
    \path (unbddbranch) -- (unbddfathom);
    \path (bddbranch) -- (bddfathom);
    \path (unbddfathom) -- (ray);
    \path (ray) -- node [midway] {\footnotesize no} (rightpt3);
    \path (rightpt3) -- (rightpt2);
    \path (rightpt2) -- (Lbarempty);
    \path (ray) -- node [midway] {\footnotesize yes} (term);
    \path (bddfathom) -- (updatelb2);
    \path (updatelb2) -- (stopbdd);
    \path (stopbdd) -- node [midway] {\footnotesize yes} (term);
    \path (stopbdd) -- node [midway] {\footnotesize no} (leftpt3);
    \path (leftpt3) -- (leftpt2);
    \path (rightpt4) -- (term);
    \path (term) -- (return);
  \end{scope}
\end{tikzpicture}
\end{center}
\caption{\label{fig. flow of whole branch-and-bound}Flow chart of branch-and-bound
procedure}
\end{figure}

The preprocessing routines are only invoked if the
the initial LP relaxation of the LPCC has a bounded optimal value;
we refer to this as the ``bounded case''.
If the initial relaxation does not have a finite optimal value then
we are in the ``unbounded case''.
For the bounded case, 
the preprocessing procedure is applied first to tighten the initial LP relaxation, then the branch-and-bound routine is invoked to solve the LPCC to optimality; for the unbounded case, we will only apply the branch-and-bound routine, which gives unbounded nodes  higher priority than  bounded nodes.
A flow diagram of the overall algorithm is given in
Figure~\ref{fig. flow of whole branch-and-bound}.
The initialization step 0
    sets the   upper bound $\bar{z}=+\infty$,
the  lower bound $\underline{z}=-\infty$,
the unbounded node list $\bar{L}=\emptyset$,
and the bounded node list $L=\emptyset$.
If the LP relaxation of the initial problem is feasible then the initial problem
is added to $L$ or $\bar{L}$ in box 1, as appropriate.
Boxes 2, 4, 6, 8, 10, 12, and 14 corresponding to the bounded
case are the subject of~\S\ref{sec: Overall Flow of Branch-and-Bound for LPCC Bounded Case},
with the unbounded case boxes 3, 5, 7, 9, and 11  explained
in~\S\ref{sec: Overall Flow of Branch-and-Bound for LPCC Unbounded Case}.
Box 13 is discussed in~\S\ref{subsec:overallBB}.

\subsection{Overall Flow of Branch-and-Bound for LPCC (Bounded Case)\label{sec: Overall Flow of Branch-and-Bound for LPCC Bounded Case} }

In the bounded case, the algorithm is quite similar to the branch-and-bound routine for a mixed integer program.
If it is determined in box 1 that the initial LP relaxation is bounded
then we implement a more detailed preprocessing
step in box 2, as discussed in~\S\ref{sec: Preprocessing Phase}.
In box 4, we apply {\em Best-Bound} to pick the next node $LPCC^i$
    from $L$ to be processed and delete $LPCC^i$ from~$L$.
    The node presolving procedure from \S\ref{sec: Node Pre-solving and warm start}
    is implemented in box 6.
The  branching strategy of
    Section \ref{subsec:branch_bdd} is used in box 8 to select  branching complementarity~$j$.
Fathoming and pruning is performed in box 10 as follows:
\begin{quote}
{\em Fathoming and pruning:} Generate two child nodes by enforcing either
    $y_j=0$ or $w_j=0$ and solve  LP relaxations. For each child node:
    \begin{enumerate}
    \item If  LP relaxation solution is feasible in LPCC with  objective $z^*$ then delete
     child node. Set $\bar{z}\leftarrow\min\{\bar{z},z^*\}$.
    \item If  LP relaxation is feasible with  objective $z^*<\bar{z}$ then set the lower bound of
    child node as $z^*$ and add  child node to~$L$.
    \item If  LP relaxation feasible with  objective $z^* \geq \bar{z}$ or infeasible then delete
     child node.
    \end{enumerate}
\end{quote}
The lower bound is updated in box 12.
The procedure is terminated in box 14 if there are no more nodes in the set~$L$
or if the gap between the upper and lower bound is sufficiently small.

\subsection{Overall Flow of Branch-and-Bound for LPCC (Unbounded Case)\label{sec: Overall Flow of Branch-and-Bound for LPCC Unbounded Case} }
The branch-and-bound routines for solving mixed integer programs in existing MIP solvers
like CPLEX usually assume the initial LP relaxation is bounded below.
Even if the initial LP relaxation is unbounded, it is still treated as bounded below by adding an objective lower bound constraint with a very large negative number ($-10^{20}$) as its lower bound. However, our branch-and-bound routine for handling the unbounded case of the LPCC is quite different. If the LP relaxation of a node is unbounded, we will treat this node as an unbounded node and add it to the unbounded node list. If the unbounded node list is non-empty, our branch-and-bound routine will always process a node in the unbounded node list first. Notice that when we find an unbounded ray that satisfies all the complementarities, we need to check whether this is a feasible ray to the LPCC. The LPCC is feasible with unbounded objective value if and only if we find an unbounded feasible ray to the LPCC.

If the set $\bar{L}$ of unbounded nodes is empty in box 3 then we return to
the bounded case in box 4, constructing an appropriate lower bound~$\underline{z}$.
In box 5, we select the node $LPCC^i$
   that is the most recently generated from $\bar{L}$
   to be processed and delete $LPCC^i$ from~$\bar{L}$.
The  branching strategy of
    Section \ref{subsec:branch_unbdd} is used in box 7 to select  branching complementarity~$j$.
Fathoming and pruning for an unbounded node is performed in box 9 as follows:
\begin{quote}
{\em Fathoming and pruning:} Generate two child nodes by enforcing either
    $y_j=0$ or $w_j=0$ and solve  LP relaxations. For each child node:
    \begin{enumerate}
    \item If  LP relaxation solution is feasible in LPCC with  objective $z^*$ then delete
     child node. Set $\bar{z}\leftarrow\min\{\bar{z},z^*\}$.
    \item If  LP relaxation is feasible with  objective $z^*<\bar{z}$ then set the lower bound of
    child node as $z^*$ and add  child node to the bounded node list~$L$.
    \item If  LP relaxation feasible with  objective $z^* \geq \bar{z}$ or infeasible then delete
    child node.
    \item If LP relaxation is unbounded and the unbounded ray is not a feasible ray to LPCC then
    add this child node to the unbounded node list~$\bar{L}$.
    \item If LP relaxation is unbounded and the piece of LPCC corresponding to that ray is feasible
    then the LPCC is unbounded.
    \end{enumerate}
\end{quote}
If an unbounded piece is found in box 11 then the algorithm can be terminated;
otherwise we loop back to box~3.

\subsection{The Complete Overall Scheme   \label{subsec:overallBB}}

A flow chart of the algorithm
is exhibited in Figure \ref{fig. flow of whole branch-and-bound}.
Each of the three possible problem states can be returned in the termination box 13.
 If an unbounded feasible ray to the LPCC is found then the LPCC is feasible
 with unbounded objective value.
    If the LPCC is not unbounded and an LPCC feasible solution is found then the LPCC attains a finite optimal   solution with optimal objective~$\bar{z}$.
    Otherwise, the problem is infeasible.

\section{Computational Results\label{sec: Computational Result for branch-and-cut algorithm}}

In this section, we will present the computational results of using our proposed branch-and-cut algorithm to solve various LPCC instances.
All procedures and algorithms are developed in the C language with the CPLEX callable library, and all LPs and convex quadratic constraint programs are solved using CPLEX 12.6.2.
We implement our algorithm through the addition of callback routines to CPLEX.
As an alternative to our approach,
CPLEX allows the modeling of complementarity constraints through the use
of {\em indicator constraints};
we compare the computational performance of our algorithm 
with that of using default CPLEX 12.6.2 to solve indicator constraint
formulations of these LPCC instances,
with our preprocessor used for both approaches.
Except for a few preliminary tests discussed in \S\ref{sec: Computational Result of branch-and-cut},
all the computational testing is performed on
a Mac Pro with 6 dual processor Intel Xeon E5 cores and 16GB of memory.
Our branch-and-cut routine uses just one thread,
while the default CPLEX 12.6.2 indicator constraint formulation can use all 12 available threads.
The relative gap for optimality is $10^{-6}$, here the relative gap is defined as $\dfrac{upperbound-lowerbound}{\mbox{max}(1,|lowerbound|)}$.
This is smaller than CPLEX's default MIP optimality tolerance
and larger than its default LP tolerance.
The tolerance of complementarity is $10^{-6}$, i.e., either $y_i$ or $w_i$ for $i=1,...,m$ should be less than $10^{-6}$ for any feasible LPCC solution.
All runtimes are reported in seconds.

We used three sets of test instances.
The first set consists of 60 LPCC instances with $n=2$ and between
100 and 200 complementarities. The generation scheme for these problems and computational results
can be found in Appendix~\ref{appendix:test instances}, with the results discussed in sections
\ref{sec: Computational Result of the Feasibility Recovery Process }
and~\ref{sec: Computational Result of branch-and-cut}.
The second set of test instances are LPCC formulations of bilevel programs,
where the lower level problem is a convex quadratic program;
the formulation and results are presented in Section~\ref{subsec:bilevel},
with more extensive results in Appendix~\ref{app:bilevel}.
The final set of results in Section~\ref{subsec:inverseQP} are for
inverse quadratic programming problems,
with detailed results in Appendix~\ref{app:inverse}.

Source code and test instances can be found online at
\url{https://github.com/mitchjrpi/LPCCbnc}
Also included with the source code is a Makefile.
A user needs to have access to CPLEX in order to be able to compile the code.
Generators for the bilevel and inverse QP problems can be found on the website;
the generator uses AMPL to construct the instances.

\subsection{Computational Results of the Feasibility Recovery Process \label{sec: Computational Result of the Feasibility Recovery Process } }
We will first apply the local search feasibility recovery process (procedure \ref{proc: Local search feasibility recovery process}); if this procedure successfully recovers a feasible solution, then the refinement procedure (procedure \ref{proc: Refined local search feasibility recovery process}) will be applied to refine that feasible solution.
We set the depth parameter as $5$ and breadth parameter as $m$, i.e. the number of complementarities, in procedure \ref{proc: Local search feasibility recovery process}.
Table \ref{table.computational Result of mean Feasibility Recovery} summarizes the feasibility recovery result of the 60 LPCC instances.
The computational results show that our proposed feasibility recovery procedures can successfully recover a feasible solution for all of the 60 LPCC instances with very good quality. For most instances, the recovered feasible solution is in fact an optimal solution. Note that as $m$ increases, the feasibility recovery processing time increases as well. Therefore in practice, as a preprocessing procedure, we need to control the depth and breadth parameters in procedure \ref{proc: Refined local search feasibility recovery process}
to reduce the time spent on the feasibility recovery procedure.

\begin{table}
\begin{center}
{
\begin{tabular}{|c|c|c|c|}
\hline
$m$ & $rankM$  & Average gap & Optimal found out of 10 \\
\hline
100 & 30 & 0.09\% & 5 \\
100 & 60 & 0.22\% & 2 \\
150 & 30 & 0.0 \% & 10\\
150 & 100 & 0.06\% & 3 \\
200 & 30 & 0.0 \% & 10 \\
200 & 120 & 0.07\% & 2 \\
\hline
\end{tabular}
}
\end{center}
\caption{\label{table.computational Result of mean Feasibility Recovery}Average Computational Results of Feasibility Recovery with $n=2$, $k=20$.
The column ``Average gap" is calculated as $\dfrac{LB_{recovered}-LPCC_{opt}}{LPCC_{opt}}$.
Detailed results can be found in Table~\ref{table:data60LPCC}
in Appendix~\ref{appendix:test instances}.
}
\end{table}

\subsection{Computational Results of Branch-and-Cut Algorithm \label{sec: Computational Result of branch-and-cut} }

In this section, we will show the computational results of using our proposed branch-and-cut algorithm to solve the 60 LPCC instances with finite global optimal values
from Appendix~\ref{appendix:test instances}.

We conducted preliminary experiments with 4 different weight settings
of increasing sophistication
to choose a score function (\ref{eq: score function of branching complementarity selection}):
\begin{enumerate}
\item[$R_1$:]$\omega^{VL}=1$, $\omega^{ED}=0$, $\omega^{PC}=0$ and $\omega^{SL}=0$;
\item[$R_2$:]$\omega^{VL}=1$, $\omega^{ED}=0.5$, $\omega^{PC}=0$ and $\omega^{SL}=0.5$;
\item[$R_3$:]$\omega^{VL}=1$, $\omega^{ED}=0.5$, $\omega^{PC}=0.25$ and $\omega^{SL}=0.5$;
\item[$R_4$:]$\omega^{VL}=1$, $\omega^{ED}=0.5$, $\omega^{PC}=0.25$ and $\omega^{SL}=0.5$ and apply strong branching rule to the node whose tree depth is less or equal to 7.
\end{enumerate}
\begin{table}
\begin{center}
{
\begin{tabular}{|c|r|r|r|r|r|}
\hline
$m$& $Time_{R_1}$ & $Time_{R_2}$ & $Time_{R_3}$ & $Time_{R_4}$ & $Time_{CPLEX}$\\ 
$ $& $(sec)$ & (sec) & (sec) & (sec) & (sec)\\ 
\hline
100&19.185&18.790&18.805&18.917&38.021\\
150&75.213&73.623&76.071&74.285&1688.160\\
200&308.293&296.663&287.232&291.057&5043.017\\
\hline
\end{tabular}
}
\end{center}
\caption{\label{table.compare solving time for 60 instances}Comparison of geometric means of solving time, using our four different branching rules and using default CPLEX}
\end{table}
\begin{table}
\begin{center}
{
\begin{tabular}{|c|r|r|r|r|r|}
\hline
$m$& $Node_{R_1}$ & $Node_{R_2}$ & $Node_{R_3}$ & $Node_{R_4}$ & $Node_{CPLEX}$\\ 
\hline
100&213&215&212&196&39114\\
150&831&863&848&757&1078311\\
200&3408&3301&3161&2837&1494577\\
\hline
\end{tabular}
}
\end{center}
\caption{\label{table.compare solving node for 60 instances}Comparison of geometric means of number of nodes in branch-and-cut tree,
using our four different branching rules and using default CPLEX}
\end{table}
These results were obtained using CPLEX 11.4 using a single core of AMD Phenom II X4 955 CPU @ 3.2GHZ, 4GB memory
and are contained in Tables \ref{table.compare solving time for 60 instances}
and~\ref{table.compare solving node for 60 instances}.
All four rules required far fewer nodes than default CPLEX.
Based on these results, $R_4$ is the best branching rule in terms of the number of  nodes.
Since in terms of solving time, these 4 routines are quite close, we chose $R_4$ as our default branch-and-bound routine.

\begin{figure}
\begin{center}
\begin{tabular}{cc}
\begin{tikzpicture}[scale=0.45]

\draw[thick] (0,12) -- (12,12) -- (12,0);
\draw[thick] (0,0) -- (12,12);
\draw[thick] (0,0) -- (0,12) 
  node[midway, sloped,above] {our algorithm};
\draw[thick,->] (0,0) -- (12,0) node[midway, below] {CPLEX};

\draw (0,0) node [below left] { 0 };
\draw (12,0) node [below] { 3600.0 };

\draw (0,12) node [left] { 3600.0 };

\fill (0.0,0.0) circle(0.1);
\fill (0.02,0.0) circle(0.1);
\fill (0.01,0.0) circle(0.1);
\fill (0.02,0.0) circle(0.1);
\fill (12.0,0.05) circle(0.1);
\fill (12.0,0.34) circle(0.1);
\fill (0.01,0.0) circle(0.1);
\fill (0.05,0.01) circle(0.1);
\fill (12.0,0.94) circle(0.1);
\fill (12.0,0.59) circle(0.1);
\fill (0.02,0.0) circle(0.1);
\fill (12.0,5.15) circle(0.1);
\fill (0.0,0.0) circle(0.1);
\fill (12.0,0.02) circle(0.1);
\fill (0.02,0.01) circle(0.1);
\fill (1.28,0.0) circle(0.1);
\fill (7.79,0.05) circle(0.1);
\fill (0.02,0.0) circle(0.1);
\fill (0.01,0.0) circle(0.1);
\fill (0.01,0.0) circle(0.1);
\fill (12.0,3.39) circle(0.1);
\fill (12.0,1.07) circle(0.1);
\fill (0.05,0.0) circle(0.1);
\fill (0.02,0.0) circle(0.1);
\fill (0.02,0.0) circle(0.1);
\fill (0.02,0.0) circle(0.1);
\fill (0.02,0.0) circle(0.1);
\fill (0.0,0.0) circle(0.1);
\fill (0.04,0.01) circle(0.1);
\fill (12.0,0.1) circle(0.1);
\fill (0.01,0.0) circle(0.1);
\fill (0.02,0.0) circle(0.1);
\fill (12.0,0.04) circle(0.1);
\fill (12.0,0.09) circle(0.1);
\fill (5.13,0.07) circle(0.1);
\fill (0.03,0.0) circle(0.1);
\fill (0.0,0.0) circle(0.1);
\fill (0.01,0.0) circle(0.1);
\fill (0.01,0.0) circle(0.1);
\fill (0.03,0.0) circle(0.1);
\fill (0.07,0.02) circle(0.1);
\fill (12.0,0.44) circle(0.1);
\fill (0.02,0.01) circle(0.1);
\fill (0.03,0.0) circle(0.1);
\fill (12.0,0.11) circle(0.1);
\fill (12.0,0.16) circle(0.1);
\fill (0.02,0.0) circle(0.1);
\fill (0.1,0.0) circle(0.1);
\fill (0.0,0.0) circle(0.1);
\fill (0.01,0.0) circle(0.1);
\fill (0.01,0.0) circle(0.1);
\fill (12.0,0.01) circle(0.1);
\fill (6.88,0.11) circle(0.1);
\fill (12.0,0.04) circle(0.1);
\fill (0.07,0.01) circle(0.1);
\fill (2.28,0.02) circle(0.1);
\fill (12.0,1.64) circle(0.1);
\fill (12.0,0.29) circle(0.1);
\fill (0.06,0.01) circle(0.1);
\fill (0.02,0.0) circle(0.1);

\end{tikzpicture}

&

\begin{tikzpicture}[scale=0.45]

\draw[thick] (0,12) -- (12,12) -- (12,0);
\draw[thick] (0,0) -- (12,12);
\draw[thick] (0,0) -- (0,12) 
  node[midway, sloped,above] {our algorithm};
\draw[thick,->] (0,0) -- (12,0) node[midway, below] {CPLEX};

\draw (0,0) node [below left] { 0 };
\draw (12,0) node [below] { 29.5 };

\draw (0,12) node [left] { 29.5 };

\fill (0.3,0.01) circle(0.1);
\fill (0.33,0.1) circle(0.1);
\fill (0.42,0.03) circle(0.1);
\fill (0.44,0.04) circle(0.1);
\fill (0.47,0.04) circle(0.1);
\fill (0.69,0.05) circle(0.1);
\fill (1.12,0.08) circle(0.1);
\fill (1.15,0.15) circle(0.1);
\fill (1.18,0.41) circle(0.1);
\fill (1.33,0.08) circle(0.1);
\fill (1.47,0.15) circle(0.1);
\fill (1.52,0.15) circle(0.1);
\fill (1.7,0.19) circle(0.1);
\fill (1.74,0.39) circle(0.1);
\fill (1.86,0.06) circle(0.1);
\fill (2.03,0.25) circle(0.1);
\fill (2.13,0.22) circle(0.1);
\fill (2.24,0.09) circle(0.1);
\fill (2.24,0.35) circle(0.1);
\fill (2.26,0.14) circle(0.1);
\fill (2.34,0.28) circle(0.1);
\fill (2.35,0.14) circle(0.1);
\fill (2.57,0.22) circle(0.1);
\fill (2.86,0.65) circle(0.1);
\fill (2.87,0.48) circle(0.1);
\fill (2.88,0.1) circle(0.1);
\fill (3.04,0.84) circle(0.1);
\fill (3.43,0.48) circle(0.1);
\fill (3.47,0.25) circle(0.1);
\fill (3.79,0.43) circle(0.1);
\fill (4.74,0.95) circle(0.1);
\fill (6.14,0.48) circle(0.1);
\fill (6.34,0.98) circle(0.1);
\fill (7.82,1.14) circle(0.1);
\fill (8.03,1.97) circle(0.1);
\fill (9.07,1.51) circle(0.1);
\fill (12.0,0.6) circle(0.1);

\end{tikzpicture}
\\
(a) & (b)
\end{tabular}
\end{center}
\caption{\label{fig:scatter_bin}Scatter plots for CPU time (seconds) for solution of LPCCs.
Horizontal axis is time for the default CPLEX indicator constraint solver, vertical axis is time for our branch and
cut algorithm.
Processing times excluded.
(a) All 60 instances.
(b) 37 LPCCs where default CPLEX MIP required no more than 30 seconds.
}
\end{figure}
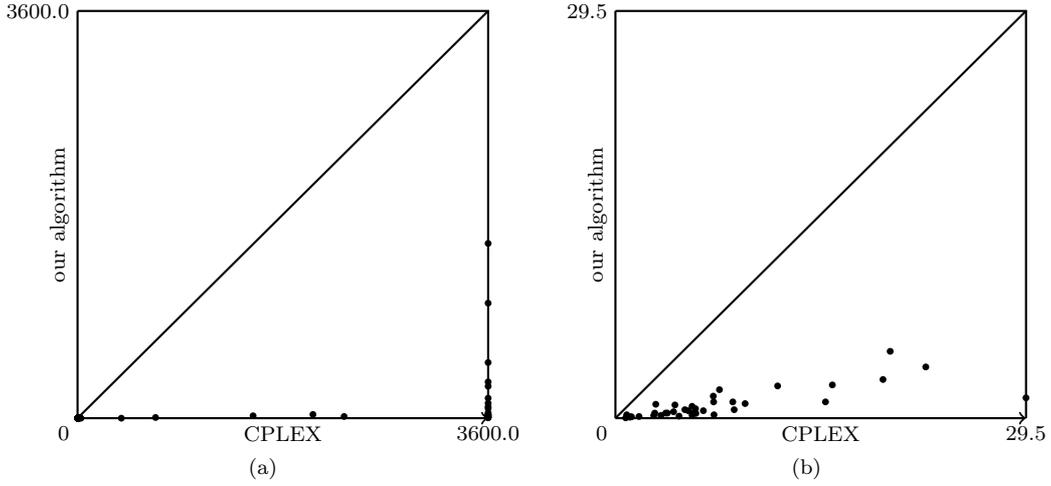

\begin{figure}
\begin{center}
\begin{tabular}{ccl}
\begin{tikzpicture}[scale=0.45]

\draw[thick,->] (0,0) -- (0,10.5) 
  node[midway, sloped,above] {$\#$ instances};

\draw[thick,->] (0,0) -- (12,0) node[midway, below] {ratio};

\draw (11,0) node [below] { 100.0 };

\draw (11,-0.05) -- (11,0.05);
\draw (0,0) node [below] { 1 };
\draw (0,0) node [left] { 0 };
\draw (0,9.84) node [left]  { 60 };

\draw[thick, blue] (0,0) -- (0.20868,0.0) -- (0.20868,0.164) -- (0.26307,0.164) -- (0.26307,0.328) -- (0.29082,0.328) -- (0.29082,0.492) -- (0.34188,0.492) -- (0.34188,0.656) -- (0.37740000000000007,0.656) -- (0.37740000000000007,0.8200000000000001) -- (0.37962,0.8200000000000001) -- (0.37962,0.984) -- (0.44621999999999995,0.984) -- (0.44621999999999995,1.1480000000000001) -- (0.54501,1.1480000000000001) -- (0.54501,1.312) -- (0.55833,1.312) -- (0.55833,1.476) -- (0.60828,1.476) -- (0.60828,1.6400000000000001) -- (0.60939,1.6400000000000001) -- (0.60939,1.804) -- (0.65268,1.804) -- (0.65268,1.968) -- (0.6782100000000001,1.968) -- (0.6782100000000001,2.132) -- (0.7437,2.132) -- (0.7437,2.2960000000000003) -- (0.78921,2.2960000000000003) -- (0.78921,2.46) -- (0.80919,2.46) -- (0.80919,2.624) -- (0.8635799999999999,2.624) -- (0.8635799999999999,2.7880000000000003) -- (0.9024300000000001,2.7880000000000003) -- (0.9024300000000001,2.952) -- (0.95571,2.952) -- (0.95571,3.116) -- (1.00566,3.116) -- (1.00566,3.2800000000000002) -- (1.0123199999999999,3.2800000000000002) -- (1.0123199999999999,3.444) -- (1.04562,3.444) -- (1.04562,3.608) -- (1.1877,3.608) -- (1.1877,3.7720000000000002) -- (1.19103,3.7720000000000002) -- (1.19103,3.936) -- (1.2654,3.936) -- (1.2654,4.1000000000000005) -- (1.31868,4.1000000000000005) -- (1.31868,4.264) -- (1.33422,4.264) -- (1.33422,4.428) -- (1.42302,4.428) -- (1.42302,4.5920000000000005) -- (1.44633,4.5920000000000005) -- (1.44633,4.756) -- (1.51626,4.756) -- (1.51626,4.92) -- (1.65279,4.92) -- (1.65279,5.0840000000000005) -- (1.6805400000000001,5.0840000000000005) -- (1.6805400000000001,5.248) -- (1.81374,5.248) -- (1.81374,5.412) -- (2.10789,5.412) -- (2.10789,5.5760000000000005) -- (2.69619,5.5760000000000005) -- (2.69619,5.74) -- (3.08913,5.74) -- (3.08913,5.904) -- (3.21123,5.904) -- (3.21123,6.0680000000000005) -- (6.94083,6.0680000000000005) -- (6.94083,6.232) -- (7.62792,6.232) -- (7.62792,6.396) -- (11.1,6.396);

\draw[thick, red, dotted] (0,0) -- (0.0,0.164) -- (0.0,0.328) -- (0.0,0.492) -- (0.0,0.656) -- (0.0,0.8200000000000001) -- (0.0,0.984) -- (0.0,1.1480000000000001) -- (0.0,1.312) -- (0.0,1.476) -- (0.0,1.6400000000000001) -- (0.0,1.804) -- (0.0,1.968) -- (0.0,2.132) -- (0.0,2.2960000000000003) -- (0.0,2.46) -- (0.0,2.624) -- (0.0,2.7880000000000003) -- (0.0,2.952) -- (0.0,3.116) -- (0.0,3.2800000000000002) -- (0.0,3.444) -- (0.0,3.608) -- (0.0,3.7720000000000002) -- (0.0,3.936) -- (0.0,4.1000000000000005) -- (0.0,4.264) -- (0.0,4.428) -- (0.0,4.5920000000000005) -- (0.0,4.756) -- (0.0,4.92) -- (0.0,5.0840000000000005) -- (0.0,5.248) -- (0.0,5.412) -- (0.0,5.5760000000000005) -- (0.0,5.74) -- (0.0,5.904) -- (0.0,6.0680000000000005) -- (0.0,6.232) -- (0.0,6.396) -- (0.0,6.5600000000000005) -- (0.0,6.724) -- (0.0,6.888) -- (0.0,7.0520000000000005) -- (0.0,7.216) -- (0.0,7.38) -- (0.0,7.5440000000000005) -- (0.0,7.708) -- (0.0,7.872) -- (0.0,8.036) -- (0.0,8.200000000000001) -- (0.0,8.364) -- (0.0,8.528) -- (0.0,8.692) -- (0.0,8.856) -- (0.0,9.02) -- (0.0,9.184000000000001) -- (0.0,9.348) -- (0.0,9.512) -- (0.0,9.676) -- (0.0,9.84) -- (11.1,9.84);

\end{tikzpicture}

&

\begin{tikzpicture}[scale=0.45]

\draw[thick,->] (0,0) -- (0,10.5)
  node[midway, sloped,above] {$\#$ instances};

\draw[thick,->] (0,0) -- (12,0) node[midway, below] {$log_{10}$(ratio)};

\draw (11,0) node [below] { 2.0 };

\draw (11,-0.05) -- (11,0.05);
\draw (0,0) node [below left] { 0 };
\draw (0,9.84) node [left]  { 60 };

\draw[thick, blue] (0,0) -- (2.5266586826757695,0.0) -- (2.5266586826757695,0.164) -- (2.9019644547923624,0.164) -- (2.9019644547923624,0.328) -- (3.0728971379324115,0.328) -- (3.0728971379324115,0.492) -- (3.3586308969943395,0.492) -- (3.3586308969943395,0.656) -- (3.5389897206740306,0.656) -- (3.5389897206740306,0.8200000000000001) -- (3.549822481420005,0.8200000000000001) -- (3.549822481420005,0.984) -- (3.8538704442976064,0.984) -- (3.8538704442976064,1.1480000000000001) -- (4.243731144846905,1.1480000000000001) -- (4.243731144846905,1.312) -- (4.291745216770832,1.312) -- (4.291745216770832,1.476) -- (4.4636625322882635,1.476) -- (4.4636625322882635,1.6400000000000001) -- (4.467345832402031,1.6400000000000001) -- (4.467345832402031,1.804) -- (4.606736410295312,1.804) -- (4.606736410295312,1.968) -- (4.685282804013715,1.968) -- (4.685282804013715,2.132) -- (4.87569898844865,2.132) -- (4.87569898844865,2.2960000000000003) -- (4.999614698161357,2.2960000000000003) -- (4.999614698161357,2.46) -- (5.052049918026504,2.46) -- (5.052049918026504,2.624) -- (5.189219837483564,2.624) -- (5.189219837483564,2.7880000000000003) -- (5.282589276438644,2.7880000000000003) -- (5.282589276438644,2.952) -- (5.404978632176999,2.952) -- (5.404978632176999,3.116) -- (5.514288893959497,3.116) -- (5.514288893959497,3.2800000000000002) -- (5.528492818770792,3.2800000000000002) -- (5.528492818770792,3.444) -- (5.598272454299281,3.444) -- (5.598272454299281,3.608) -- (5.875022239603889,3.608) -- (5.875022239603889,3.7720000000000002) -- (5.881139066635411,3.7720000000000002) -- (5.881139066635411,3.936) -- (6.013819268392293,3.936) -- (6.013819268392293,4.1000000000000005) -- (6.104537246630863,4.1000000000000005) -- (6.104537246630863,4.264) -- (6.130360413276953,4.264) -- (6.130360413276953,4.428) -- (6.272794236709988,4.428) -- (6.272794236709988,4.5920000000000005) -- (6.30881719065598,4.5920000000000005) -- (6.30881719065598,4.756) -- (6.413736836678101,4.756) -- (6.413736836678101,4.92) -- (6.606181434640588,4.92) -- (6.606181434640588,5.0840000000000005) -- (6.643469417123284,5.0840000000000005) -- (6.643469417123284,5.248) -- (6.814770012271053,5.248) -- (6.814770012271053,5.412) -- (7.154470367649643,5.412) -- (7.154470367649643,5.5760000000000005) -- (7.716218561394227,5.5760000000000005) -- (7.716218561394227,5.74) -- (8.029145533135143,5.74) -- (8.029145533135143,5.904) -- (8.118586942621052,5.904) -- (8.118586942621052,6.0680000000000005) -- (9.916383704343941,6.0680000000000005) -- (9.916383704343941,6.232) -- (10.138465581408756,6.232) -- (10.138465581408756,6.396) -- (11.1,6.396);

\draw[thick, red, dotted] (0,0) -- (0.0,0.164) -- (0.0,0.328) -- (0.0,0.492) -- (0.0,0.656) -- (0.0,0.8200000000000001) -- (0.0,0.984) -- (0.0,1.1480000000000001) -- (0.0,1.312) -- (0.0,1.476) -- (0.0,1.6400000000000001) -- (0.0,1.804) -- (0.0,1.968) -- (0.0,2.132) -- (0.0,2.2960000000000003) -- (0.0,2.46) -- (0.0,2.624) -- (0.0,2.7880000000000003) -- (0.0,2.952) -- (0.0,3.116) -- (0.0,3.2800000000000002) -- (0.0,3.444) -- (0.0,3.608) -- (0.0,3.7720000000000002) -- (0.0,3.936) -- (0.0,4.1000000000000005) -- (0.0,4.264) -- (0.0,4.428) -- (0.0,4.5920000000000005) -- (0.0,4.756) -- (0.0,4.92) -- (0.0,5.0840000000000005) -- (0.0,5.248) -- (0.0,5.412) -- (0.0,5.5760000000000005) -- (0.0,5.74) -- (0.0,5.904) -- (0.0,6.0680000000000005) -- (0.0,6.232) -- (0.0,6.396) -- (0.0,6.5600000000000005) -- (0.0,6.724) -- (0.0,6.888) -- (0.0,7.0520000000000005) -- (0.0,7.216) -- (0.0,7.38) -- (0.0,7.5440000000000005) -- (0.0,7.708) -- (0.0,7.872) -- (0.0,8.036) -- (0.0,8.200000000000001) -- (0.0,8.364) -- (0.0,8.528) -- (0.0,8.692) -- (0.0,8.856) -- (0.0,9.02) -- (0.0,9.184000000000001) -- (0.0,9.348) -- (0.0,9.512) -- (0.0,9.676) -- (0.0,9.84) -- (11.1,9.84);

\end{tikzpicture}

&

\begin{tikzpicture}[scale=0.3]
\draw[thick, red, dotted] (0,9) -- (0.5,9);
\draw (0.8,9) node [right] {branch-and-cut};
\draw[thick, blue] (0,7) -- (0.5,7);
\draw (0.8,7) node [right] {CPLEX indicator};
\draw (0.82,5.8) node [right] {constraint};
\draw (0,0) node [right] {$\quad$};
\end{tikzpicture}

\\
(a) & (b)
\end{tabular}
\end{center}
\caption{\label{fig:perfprofBin_macpro}Performance profile for CPU time (seconds) for solution of 60 LPCCs
(preprocessing time excluded).
Vertical axis is the number of instances.
Horizontal axis is ratio of time required by the given algorithm to the
time required by the better algorithm.
(a) Linear scale.  (b) Log scale.}
\end{figure}
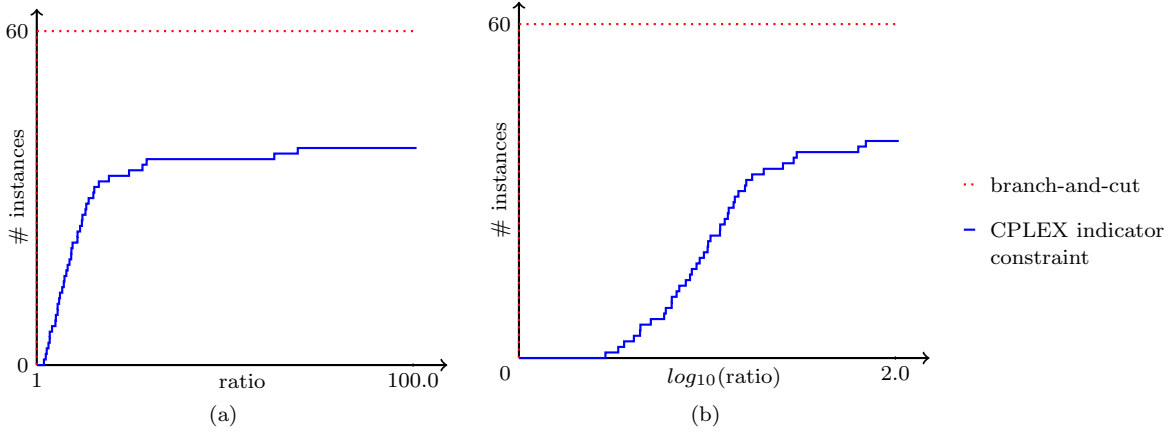

All remaining results in the paper were obtained using CPLEX 12.6.2
with detailed results contained in Table~\ref{table:results60LPCC}
in Appendix~\ref{appendix:test instances}.
A scatter plot of the CPU time for solving the instances is given in Figure~\ref{fig:scatter_bin}.
Performance profiles~\cite{dolan2} are given in Figure~\ref{fig:perfprofBin_macpro}.
The preprocessing times have been excluded from these plots.
All the LPCC instances can be solved by our algorithm within thirty minutes,
with 90\% of them (54/60)  solved within 150 seconds.
Each instance requires considerably less processing time with our algorithm than with default
CPLEX.
Notice that default CPLEX is only able to solve 42 of the 60 instances within 3600 seconds.
In particular, it is unable to solve 11 of our 20 LPCC instances when $m=200$ within
this time limit.

The determination of a valid disjunctive cut or bound cut requires
the solution of a linear programming problem.
The parameter choices given
in~\S\ref{sec: Cutting Planes Generation and Selection}
result in $0.3m$ disjunctive cuts, approximately $5m$ simple cuts,
and 15 bound cuts for each instance.
We also experimented with not adding cutting planes in the preprocessor,
in which case
both codes performed slightly worse for the larger instances
(a difference of perhaps 10\% in average runtime for our branch-and-cut code).

\subsection{Bilevel Test Problems   \label{subsec:bilevel}}

We further tested our algorithm on bilevel problems of the form
\begin{equation}  \label{eq:bilevel}
\begin{array}{lrcrcl}
\min_{x,v} & c^Tx & + & d^Tv  \\
\mbox{subject to} & Ax & + & Bv & \geq & b  \\
& 0 & \leq & v & \leq & u \\
& x & \in & \multicolumn{3}{l}{\mbox{argmin}_x\{ \frac{1}{2} x^TQx \, + \, v^Tx \, : \, Hx \, \geq \, g, \, x \geq 0\} }
\end{array}
\end{equation}
where $Q$ is positive semidefinite.
The variables $v$ are first stage variables,
with the second stage variables $x$ chosen to optimize a convex quadratic subproblem
that depends on~$v$.
Both sets of variables appear in the linear objective.
In addition, the first and second stage variables must satisfy the linking
constraint $Ax+Bv \geq b$.
By introducing KKT multipliers $y$ and $\lambda$ for the constraints in the subproblem,
we can model this problem equivalently as the LPCC
\begin{displaymath}
\begin{array}{lrcrclcl}
\min_{x,v,y,\lambda,w} & c^Tx & + & d^Tv  \\
\mbox{subject to} & Ax & + & Bv &&& \geq & b  \\
& Qx & + & v & - & H^Ty \, - \, \lambda & = & 0  \\
& 0 & \leq & v &&& \leq & u \\
& 0 & \leq & \lambda & \perp & x & \geq & 0  \\
& 0 & \leq & y & \perp & w \, := \, Hx \, - \, g & \geq & 0,
\end{array}
\end{displaymath}
a problem equivalent to one in our standard form~(\ref{eq:general LPCC}).
The relationship between the dimensions in (\ref{eq:general LPCC}) and the dimensions of
the variables and constraints in (\ref{eq:bilevel}) is as follows:
\begin{center}
\begin{tabular}{|l|l|}
\multicolumn{2}{c}{Dimensions}  \\  \hline
(\ref{eq:general LPCC})  &  (\ref{eq:bilevel})  \\  \hline
$m$ & dimension($g$) $+$ dimension($v$)  \\
$n$ & $2$ $\times$ dimension($v$)  \\
$k$ & dimension($b$) $+$ $3$ $\times$ dimension($v$)  \\ \hline
\end{tabular}
\end{center}
Thus, the number of complementarity constraints is equal to the sum of the
dimensions of $g$ and~$v$.

In our experiments,
all parameters in $b$, $c$, $d$, $g$, $A$, $B$, and $H$ were uniformly
generated in the interval $(0,1)$.
The matrix $Q$ was equal to the matrix product $LL^T$,
where the number of columns in $L$ is equal to the required rank of~$Q$
and each entry in $L$ is chosen uniformly from the interval $(-1,1)$.
Each entry of $u$ was equal to~1.
Repeated problem dimensions in the table correspond to different
instances.
The dimension of $g$ varied from 50 to 200,
the dimension of $v$ and $x$ varied from 50 to 100,
the number of complementarity constraints varied from 100 to 250,
the dimension of $b$ varied from 25 to 100,
and the rank of $Q$ varied between 0.5 of the dimension of $v$ and the dimension of~$v$. 
Problem data for the 90 bilevel test instances can be found in
Tables~\ref{table:vals_bil_50} and~\ref{table:vals_bil_larger}.

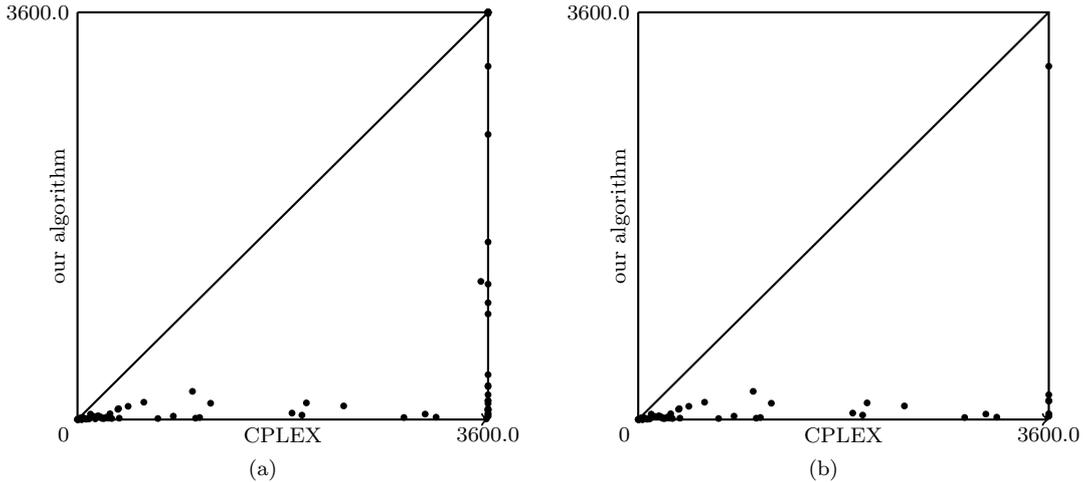
\begin{figure}
\begin{center}
\begin{tabular}{cc}
\begin{tikzpicture}[scale=0.45]

\draw[thick] (0,12) -- (12,12) -- (12,0);
\draw[thick] (0,0) -- (12,12);
\draw[thick] (0,0) -- (0,12) 
  node[midway, sloped,above] {our algorithm};
\draw[thick,->] (0,0) -- (12,0) node[midway, below] {CPLEX};

\draw (0,0) node [below left] { 0 };
\draw (12,0) node [below] { 3600.0 };

\draw (0,12) node [left] { 3600.0 };

\fill (12.0,12.0) circle(0.1);
\fill (0.52,0.02) circle(0.1);
\fill (12.0,12.0) circle(0.1);
\fill (12.0,12.0) circle(0.1);
\fill (12.0,12.0) circle(0.1);
\fill (7.78,0.4) circle(0.1);
\fill (0.49,0.09) circle(0.1);
\fill (12.0,0.73) circle(0.1);
\fill (10.16,0.16) circle(0.1);
\fill (12.0,10.41) circle(0.1);
\fill (0.2,0.03) circle(0.1);
\fill (0.91,0.04) circle(0.1);
\fill (0.12,0.0) circle(0.1);
\fill (0.38,0.16) circle(0.1);
\fill (0.13,0.01) circle(0.1);
\fill (12.0,0.11) circle(0.1);
\fill (2.35,0.03) circle(0.1);
\fill (12.0,0.54) circle(0.1);
\fill (0.01,0.0) circle(0.1);
\fill (0.21,0.02) circle(0.1);
\fill (0.08,0.02) circle(0.1);
\fill (12.0,0.1) circle(0.1);
\fill (0.28,0.02) circle(0.1);
\fill (1.2,0.32) circle(0.1);
\fill (0.0,0.0) circle(0.1);
\fill (3.36,0.83) circle(0.1);
\fill (0.02,0.01) circle(0.1);
\fill (0.95,0.17) circle(0.1);
\fill (0.73,0.03) circle(0.1);
\fill (0.09,0.05) circle(0.1);
\fill (0.0,0.0) circle(0.1);
\fill (0.01,0.01) circle(0.1);
\fill (9.54,0.06) circle(0.1);
\fill (6.27,0.19) circle(0.1);
\fill (0.03,0.01) circle(0.1);
\fill (0.15,0.03) circle(0.1);
\fill (12.0,1.32) circle(0.1);
\fill (12.0,12.0) circle(0.1);
\fill (12.0,0.97) circle(0.1);
\fill (12.0,12.0) circle(0.1);
\fill (12.0,0.17) circle(0.1);
\fill (12.0,12.0) circle(0.1);
\fill (12.0,0.3) circle(0.1);
\fill (12.0,11.96) circle(0.1);
\fill (12.0,8.4) circle(0.1);
\fill (12.0,12.0) circle(0.1);
\fill (12.0,0.31) circle(0.1);
\fill (12.0,0.26) circle(0.1);
\fill (12.0,3.44) circle(0.1);
\fill (12.0,12.0) circle(0.1);
\fill (0.17,0.02) circle(0.1);
\fill (6.69,0.49) circle(0.1);
\fill (0.03,0.0) circle(0.1);
\fill (0.15,0.06) circle(0.1);
\fill (1.0,0.03) circle(0.1);
\fill (1.18,0.3) circle(0.1);
\fill (0.78,0.03) circle(0.1);
\fill (12.0,0.17) circle(0.1);
\fill (1.22,0.04) circle(0.1);
\fill (2.8,0.1) circle(0.1);
\fill (3.57,0.06) circle(0.1);
\fill (0.85,0.07) circle(0.1);
\fill (1.48,0.39) circle(0.1);
\fill (12.0,0.56) circle(0.1);
\fill (0.1,0.01) circle(0.1);
\fill (0.24,0.03) circle(0.1);
\fill (0.13,0.01) circle(0.1);
\fill (0.7,0.05) circle(0.1);
\fill (0.72,0.05) circle(0.1);
\fill (1.94,0.51) circle(0.1);
\fill (6.56,0.13) circle(0.1);
\fill (3.89,0.48) circle(0.1);
\fill (0.14,0.01) circle(0.1);
\fill (0.31,0.03) circle(0.1);
\fill (3.45,0.04) circle(0.1);
\fill (0.34,0.02) circle(0.1);
\fill (10.48,0.07) circle(0.1);
\fill (0.6,0.11) circle(0.1);
\fill (0.09,0.02) circle(0.1);
\fill (0.65,0.09) circle(0.1);
\fill (0.26,0.01) circle(0.1);
\fill (0.44,0.08) circle(0.1);
\fill (12.0,3.11) circle(0.1);
\fill (12.0,1.0) circle(0.1);
\fill (12.0,3.99) circle(0.1);
\fill (11.95,0.02) circle(0.1);
\fill (12.0,0.47) circle(0.1);
\fill (12.0,5.23) circle(0.1);
\fill (11.79,4.07) circle(0.1);
\fill (12.0,0.29) circle(0.1);

\end{tikzpicture}

&

\begin{tikzpicture}[scale=0.45]

\draw[thick] (0,12) -- (12,12) -- (12,0);
\draw[thick] (0,0) -- (12,12);
\draw[thick] (0,0) -- (0,12) 
  node[midway, sloped,above] {our algorithm};
\draw[thick,->] (0,0) -- (12,0) node[midway, below] {CPLEX};

\draw (0,0) node [below left] { 0 };
\draw (12,0) node [below] { 3600.0 };

\draw (0,12) node [left] { 3600.0 };

\fill (0.17,0.02) circle(0.1);
\fill (6.69,0.49) circle(0.1);
\fill (0.03,0.0) circle(0.1);
\fill (0.15,0.06) circle(0.1);
\fill (1.0,0.03) circle(0.1);
\fill (1.18,0.3) circle(0.1);
\fill (0.78,0.03) circle(0.1);
\fill (12.0,0.17) circle(0.1);
\fill (1.22,0.04) circle(0.1);
\fill (2.8,0.1) circle(0.1);
\fill (3.57,0.06) circle(0.1);
\fill (0.85,0.07) circle(0.1);
\fill (1.48,0.39) circle(0.1);
\fill (12.0,0.56) circle(0.1);
\fill (0.1,0.01) circle(0.1);
\fill (0.24,0.03) circle(0.1);
\fill (0.13,0.01) circle(0.1);
\fill (0.7,0.05) circle(0.1);
\fill (0.72,0.05) circle(0.1);
\fill (1.94,0.51) circle(0.1);
\fill (6.56,0.13) circle(0.1);
\fill (3.89,0.48) circle(0.1);
\fill (0.14,0.01) circle(0.1);
\fill (0.31,0.03) circle(0.1);
\fill (3.45,0.04) circle(0.1);
\fill (0.34,0.02) circle(0.1);
\fill (10.48,0.07) circle(0.1);
\fill (0.6,0.11) circle(0.1);
\fill (0.09,0.02) circle(0.1);
\fill (0.65,0.09) circle(0.1);
\fill (0.26,0.01) circle(0.1);
\fill (0.44,0.08) circle(0.1);
\fill (7.78,0.4) circle(0.1);
\fill (0.49,0.09) circle(0.1);
\fill (12.0,0.73) circle(0.1);
\fill (10.16,0.16) circle(0.1);
\fill (12.0,10.41) circle(0.1);
\fill (0.2,0.03) circle(0.1);
\fill (0.91,0.04) circle(0.1);
\fill (0.12,0.0) circle(0.1);
\fill (0.38,0.16) circle(0.1);
\fill (0.13,0.01) circle(0.1);
\fill (12.0,0.11) circle(0.1);
\fill (2.35,0.03) circle(0.1);
\fill (12.0,0.54) circle(0.1);
\fill (0.01,0.0) circle(0.1);
\fill (0.21,0.02) circle(0.1);
\fill (0.08,0.02) circle(0.1);
\fill (12.0,0.1) circle(0.1);
\fill (0.28,0.02) circle(0.1);
\fill (1.2,0.32) circle(0.1);
\fill (0.0,0.0) circle(0.1);
\fill (3.36,0.83) circle(0.1);
\fill (0.02,0.01) circle(0.1);
\fill (0.95,0.17) circle(0.1);
\fill (0.73,0.03) circle(0.1);
\fill (0.09,0.05) circle(0.1);
\fill (0.0,0.0) circle(0.1);
\fill (0.01,0.01) circle(0.1);
\fill (9.54,0.06) circle(0.1);
\fill (6.27,0.19) circle(0.1);
\fill (0.03,0.01) circle(0.1);
\fill (0.15,0.03) circle(0.1);

\end{tikzpicture}

\\
(a) & (b)
\end{tabular}
\end{center}
\caption{\label{fig:scatter_bil}Scatter plots for CPU time (seconds) for solution of LPCCs
based on bilevel instances
Horizontal axis is time for default CPLEX indicator constraint solver,
vertical axis is time for our branch-and-cut solver.
Processing times excluded.
(a) All 90 instances.
(b) 63 instances with $n=50$.
}
\end{figure}

\begin{figure}
\begin{center}
\begin{tabular}{ccl}
\begin{tikzpicture}[scale=0.49]

\draw[thick,->] (0,0) -- (0,9.7)
  node[midway, sloped,above] {$\#$ instances};

\draw[thick,->] (0,0) -- (12,0) node[midway, below] {ratio};

\draw (11,0) node [below] { 100.0 };

\draw (11,-0.05) -- (11,0.05);
\draw (0,0) node [below] { 1 };
\draw (0,0) node [left] { 0 };
\draw (0,9.090000000000001) node [left]  { 90 };
\draw (0,8.282) node [left]  { 82 };

\draw[thick, blue] (0,0) -- (0.0666,0.0) -- (0.0666,0.101) -- (0.10989,0.101) -- (0.10989,0.202) -- (0.15539999999999998,0.202) -- (0.15539999999999998,0.30300000000000005) -- (0.16095,0.30300000000000005) -- (0.16095,0.404) -- (0.2109,0.404) -- (0.2109,0.505) -- (0.21312,0.505) -- (0.21312,0.6060000000000001) -- (0.21645000000000003,0.6060000000000001) -- (0.21645000000000003,0.7070000000000001) -- (0.30635999999999997,0.7070000000000001) -- (0.30635999999999997,0.808) -- (0.30747,0.808) -- (0.30747,0.909) -- (0.31301999999999996,0.909) -- (0.31301999999999996,1.01) -- (0.31968,1.01) -- (0.31968,1.111) -- (0.33965999999999996,1.111) -- (0.33965999999999996,1.2120000000000002) -- (0.34965,1.2120000000000002) -- (0.34965,1.3130000000000002) -- (0.36185999999999996,1.3130000000000002) -- (0.36185999999999996,1.4140000000000001) -- (0.38516999999999996,1.4140000000000001) -- (0.38516999999999996,1.5150000000000001) -- (0.49284000000000006,1.5150000000000001) -- (0.49284000000000006,1.616) -- (0.5061599999999999,1.616) -- (0.5061599999999999,1.717) -- (0.50949,1.717) -- (0.50949,1.818) -- (0.5372399999999999,1.818) -- (0.5372399999999999,1.919) -- (0.6260399999999999,1.919) -- (0.6260399999999999,2.02) -- (0.6304799999999999,2.02) -- (0.6304799999999999,2.121) -- (0.71484,2.121) -- (0.71484,2.222) -- (0.76257,2.222) -- (0.76257,2.323) -- (0.7914300000000001,2.323) -- (0.7914300000000001,2.4240000000000004) -- (0.8080799999999999,2.4240000000000004) -- (0.8080799999999999,2.5250000000000004) -- (0.9024300000000001,2.5250000000000004) -- (0.9024300000000001,2.6260000000000003) -- (0.90687,2.6260000000000003) -- (0.90687,2.7270000000000003) -- (0.96459,2.7270000000000003) -- (0.96459,2.8280000000000003) -- (1.03785,2.8280000000000003) -- (1.03785,2.9290000000000003) -- (1.08003,2.9290000000000003) -- (1.08003,3.0300000000000002) -- (1.10667,3.0300000000000002) -- (1.10667,3.1310000000000002) -- (1.2731700000000001,3.1310000000000002) -- (1.2731700000000001,3.232) -- (1.37085,3.232) -- (1.37085,3.333) -- (1.4019300000000001,3.333) -- (1.4019300000000001,3.434) -- (1.43856,3.434) -- (1.43856,3.535) -- (1.5695400000000002,3.535) -- (1.5695400000000002,3.636) -- (1.6805400000000001,3.636) -- (1.6805400000000001,3.737) -- (2.03796,3.737) -- (2.03796,3.838) -- (2.14119,3.838) -- (2.14119,3.939) -- (2.22777,3.939) -- (2.22777,4.04) -- (2.43201,4.04) -- (2.43201,4.141) -- (2.49639,4.141) -- (2.49639,4.242) -- (2.50638,4.242) -- (2.50638,4.343) -- (2.59185,4.343) -- (2.59185,4.444) -- (2.71062,4.444) -- (2.71062,4.545) -- (2.9992199999999998,4.545) -- (2.9992199999999998,4.646) -- (3.29448,4.646) -- (3.29448,4.747) -- (3.4543199999999996,4.747) -- (3.4543199999999996,4.848000000000001) -- (3.5852999999999997,4.848000000000001) -- (3.5852999999999997,4.949000000000001) -- (3.9704700000000006,4.949000000000001) -- (3.9704700000000006,5.050000000000001) -- (5.55666,5.050000000000001) -- (5.55666,5.151000000000001) -- (6.68109,5.151000000000001) -- (6.68109,5.252000000000001) -- (6.8376,5.252000000000001) -- (6.8376,5.353000000000001) -- (7.0773600000000005,5.353000000000001) -- (7.0773600000000005,5.454000000000001) -- (8.51703,5.454000000000001) -- (8.51703,5.555000000000001) -- (9.16305,5.555000000000001) -- (9.16305,5.656000000000001) -- (11.1,5.656000000000001);

\draw[thick, red, dotted] (0,0) -- (0.0,0.101) -- (0.0,0.202) -- (0.0,0.30300000000000005) -- (0.0,0.404) -- (0.0,0.505) -- (0.0,0.6060000000000001) -- (0.0,0.7070000000000001) -- (0.0,0.808) -- (0.0,0.909) -- (0.0,1.01) -- (0.0,1.111) -- (0.0,1.2120000000000002) -- (0.0,1.3130000000000002) -- (0.0,1.4140000000000001) -- (0.0,1.5150000000000001) -- (0.0,1.616) -- (0.0,1.717) -- (0.0,1.818) -- (0.0,1.919) -- (0.0,2.02) -- (0.0,2.121) -- (0.0,2.222) -- (0.0,2.323) -- (0.0,2.4240000000000004) -- (0.0,2.5250000000000004) -- (0.0,2.6260000000000003) -- (0.0,2.7270000000000003) -- (0.0,2.8280000000000003) -- (0.0,2.9290000000000003) -- (0.0,3.0300000000000002) -- (0.0,3.1310000000000002) -- (0.0,3.232) -- (0.0,3.333) -- (0.0,3.434) -- (0.0,3.535) -- (0.0,3.636) -- (0.0,3.737) -- (0.0,3.838) -- (0.0,3.939) -- (0.0,4.04) -- (0.0,4.141) -- (0.0,4.242) -- (0.0,4.343) -- (0.0,4.444) -- (0.0,4.545) -- (0.0,4.646) -- (0.0,4.747) -- (0.0,4.848000000000001) -- (0.0,4.949000000000001) -- (0.0,5.050000000000001) -- (0.0,5.151000000000001) -- (0.0,5.252000000000001) -- (0.0,5.353000000000001) -- (0.0,5.454000000000001) -- (0.0,5.555000000000001) -- (0.0,5.656000000000001) -- (0.0,5.757000000000001) -- (0.0,5.8580000000000005) -- (0.0,5.9590000000000005) -- (0.0,6.0600000000000005) -- (0.0,6.1610000000000005) -- (0.0,6.2620000000000005) -- (0.0,6.363) -- (0.0,6.464) -- (0.0,6.565) -- (0.0,6.666) -- (0.0,6.767) -- (0.0,6.868) -- (0.0,6.969) -- (0.0,7.07) -- (0.0,7.171) -- (0.0,7.272) -- (0.0,7.373) -- (0.0,7.474) -- (0.0,7.575) -- (0.0,7.676) -- (0.0,7.777) -- (0.0,7.878) -- (0.0,7.979) -- (0.0,8.08) -- (0.0,8.181000000000001) -- (11.1,8.181000000000001);

\end{tikzpicture}

&

\begin{tikzpicture}[scale=0.45]

\draw[thick,->] (0,0) -- (0,10.5)
  node[midway, sloped,above] {$\#$ instances};

\draw[thick,->] (0,0) -- (12,0) node[midway, below] {ratio};

\draw (11,0) node [below] { 100.0 };

\draw (11,-0.05) -- (11,0.05);
\draw (0,0) node [below] { 1 };
\draw (0,0) node [left] { 0 };
\draw (0,9.828) node [left]  { 63 };

\draw[thick, blue] (0,0) -- (0.0666,0.0) -- (0.0666,0.156) -- (0.10989,0.156) -- (0.10989,0.312) -- (0.15539999999999998,0.312) -- (0.15539999999999998,0.46799999999999997) -- (0.16095,0.46799999999999997) -- (0.16095,0.624) -- (0.21312,0.624) -- (0.21312,0.78) -- (0.21645000000000003,0.78) -- (0.21645000000000003,0.9359999999999999) -- (0.30635999999999997,0.9359999999999999) -- (0.30635999999999997,1.092) -- (0.30747,1.092) -- (0.30747,1.248) -- (0.31301999999999996,1.248) -- (0.31301999999999996,1.404) -- (0.31968,1.404) -- (0.31968,1.56) -- (0.33965999999999996,1.56) -- (0.33965999999999996,1.716) -- (0.34965,1.716) -- (0.34965,1.8719999999999999) -- (0.36185999999999996,1.8719999999999999) -- (0.36185999999999996,2.028) -- (0.38516999999999996,2.028) -- (0.38516999999999996,2.184) -- (0.49284000000000006,2.184) -- (0.49284000000000006,2.34) -- (0.5061599999999999,2.34) -- (0.5061599999999999,2.496) -- (0.50949,2.496) -- (0.50949,2.652) -- (0.5372399999999999,2.652) -- (0.5372399999999999,2.808) -- (0.6260399999999999,2.808) -- (0.6260399999999999,2.964) -- (0.6304799999999999,2.964) -- (0.6304799999999999,3.12) -- (0.71484,3.12) -- (0.71484,3.276) -- (0.76257,3.276) -- (0.76257,3.432) -- (0.7914300000000001,3.432) -- (0.7914300000000001,3.588) -- (0.8080799999999999,3.588) -- (0.8080799999999999,3.7439999999999998) -- (0.9024300000000001,3.7439999999999998) -- (0.9024300000000001,3.9) -- (0.90687,3.9) -- (0.90687,4.056) -- (0.96459,4.056) -- (0.96459,4.212) -- (1.03785,4.212) -- (1.03785,4.368) -- (1.08003,4.368) -- (1.08003,4.524) -- (1.10667,4.524) -- (1.10667,4.68) -- (1.2731700000000001,4.68) -- (1.2731700000000001,4.836) -- (1.37085,4.836) -- (1.37085,4.992) -- (1.4019300000000001,4.992) -- (1.4019300000000001,5.148) -- (1.43856,5.148) -- (1.43856,5.304) -- (1.5695400000000002,5.304) -- (1.5695400000000002,5.46) -- (1.6805400000000001,5.46) -- (1.6805400000000001,5.616) -- (2.03796,5.616) -- (2.03796,5.772) -- (2.14119,5.772) -- (2.14119,5.928) -- (2.22777,5.928) -- (2.22777,6.084) -- (2.49639,6.084) -- (2.49639,6.24) -- (2.50638,6.24) -- (2.50638,6.396) -- (2.59185,6.396) -- (2.59185,6.552) -- (2.71062,6.552) -- (2.71062,6.708) -- (2.9992199999999998,6.708) -- (2.9992199999999998,6.864) -- (3.29448,6.864) -- (3.29448,7.02) -- (3.4543199999999996,7.02) -- (3.4543199999999996,7.176) -- (3.5852999999999997,7.176) -- (3.5852999999999997,7.332) -- (3.9704700000000006,7.332) -- (3.9704700000000006,7.4879999999999995) -- (5.55666,7.4879999999999995) -- (5.55666,7.644) -- (6.68109,7.644) -- (6.68109,7.8) -- (6.8376,7.8) -- (6.8376,7.956) -- (7.0773600000000005,7.956) -- (7.0773600000000005,8.112) -- (8.51703,8.112) -- (8.51703,8.268) -- (9.16305,8.268) -- (9.16305,8.424) -- (11.1,8.424);

\draw[thick, red, dotted] (0,0) -- (0.0,0.156) -- (0.0,0.312) -- (0.0,0.46799999999999997) -- (0.0,0.624) -- (0.0,0.78) -- (0.0,0.9359999999999999) -- (0.0,1.092) -- (0.0,1.248) -- (0.0,1.404) -- (0.0,1.56) -- (0.0,1.716) -- (0.0,1.8719999999999999) -- (0.0,2.028) -- (0.0,2.184) -- (0.0,2.34) -- (0.0,2.496) -- (0.0,2.652) -- (0.0,2.808) -- (0.0,2.964) -- (0.0,3.12) -- (0.0,3.276) -- (0.0,3.432) -- (0.0,3.588) -- (0.0,3.7439999999999998) -- (0.0,3.9) -- (0.0,4.056) -- (0.0,4.212) -- (0.0,4.368) -- (0.0,4.524) -- (0.0,4.68) -- (0.0,4.836) -- (0.0,4.992) -- (0.0,5.148) -- (0.0,5.304) -- (0.0,5.46) -- (0.0,5.616) -- (0.0,5.772) -- (0.0,5.928) -- (0.0,6.084) -- (0.0,6.24) -- (0.0,6.396) -- (0.0,6.552) -- (0.0,6.708) -- (0.0,6.864) -- (0.0,7.02) -- (0.0,7.176) -- (0.0,7.332) -- (0.0,7.4879999999999995) -- (0.0,7.644) -- (0.0,7.8) -- (0.0,7.956) -- (0.0,8.112) -- (0.0,8.268) -- (0.0,8.424) -- (0.0,8.58) -- (0.0,8.736) -- (0.0,8.892) -- (0.0,9.048) -- (0.0,9.204) -- (0.0,9.36) -- (0.0,9.516) -- (0.0,9.672) -- (0.0,9.828) -- (11.1,9.828);

\end{tikzpicture}

&

\begin{tikzpicture}[scale=0.3]
\draw[thick, red, dotted] (0,9) -- (0.5,9);
\draw (0.8,9) node [right] {branch-and-cut};
\draw[thick, blue] (0,7) -- (0.5,7);
\draw (0.8,7) node [right] {CPLEX indicator};
\draw (0.82,5.8) node [right] {constraint};
\draw (0,0) node [right] {$\quad$};
\end{tikzpicture}

\\

(a) & (b)
\end{tabular}
\end{center}
\caption{\label{fig:perf_bil}Performance profile with linear scale
for CPU time (seconds) for solution of  LPCCs
based on bilevel instances
(preprocessing time excluded).
Vertical axis is the number of instances.
Horizontal axis is ratio of time required by the given algorithm to the
time required by the better algorithm.
(a) All 90 instances.  (b) 63 instances with $n=50$.}
\end{figure}
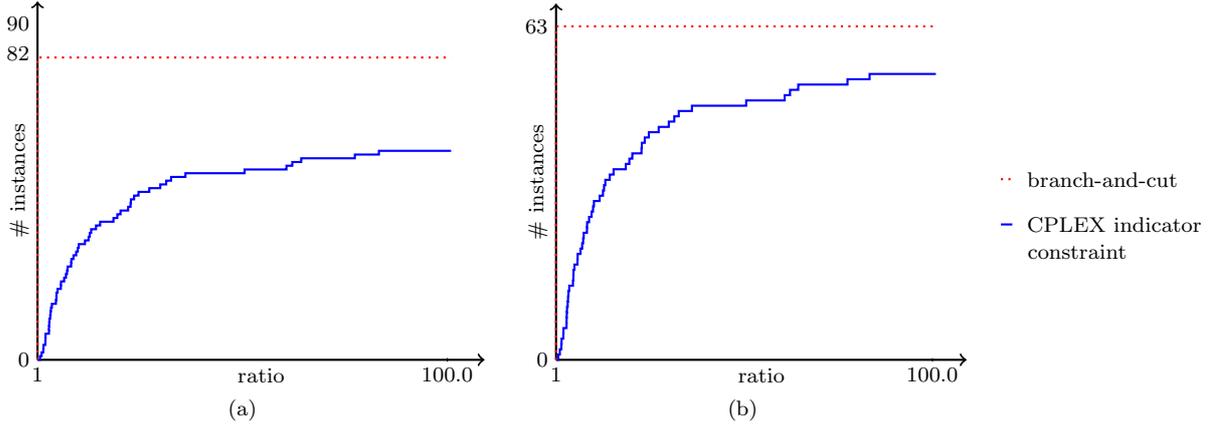

We gave each algorithm a time limit of 3600 seconds in addition to the preprocessing time.
Detailed performance data can be found in Tables \ref{table:perf_bil_50} and~\ref{table:perf_bil_larger}.
Our algorithm was able to solve all 63 instances with dimension of $v$ equal to 50,
16/24 of the instances with the dimension of $v$ equal to 75,
and 3/3 of the instances with the dimension of $v$ equal to 100.
The corresponding numbers for the default CPLEX indicator constraint code were 56/63, 2/24, and 1/3.
Our algorithm was considerably faster than default CPLEX indicator constraint code on every instance.
Further, it had a smaller final gap than default CPLEX indicator constraint code for each instance where neither code
could solve the problem.
There was no instance that could be solved by default CPLEX indicator constraint code which could not also be solved
by our algorithm.
A scatter plot of the CPU time for solving the instances (ignoring the common
preprocessing time) is given in Figure~\ref{fig:scatter_bil} and
a performance profile is in Figure~\ref{fig:perf_bil}.

The instances become more difficult as the dimensions of $v$, $b$, and $g$ increase,
as might be expected.
The instances also become more difficult as the rank of $Q$ increases.
Table~\ref{table:averages50} contains
averages of solution times over these different parameters
for the instances with the dimension of $v$ equal to 50.
\begin{table}
\begin{center}
\begin{tabular}{|r|r|r||r|r|r||r|r|r|}
\hline
\multicolumn{3}{|c||}{Dimension of $g$} & \multicolumn{3}{c||}{Dimension of $b$} &
\multicolumn{3}{c|}{Rank of $Q$}  \\ \hline
dim & time & \# instances &
dim & time & \# instances &
rank & time & \# instances \\  \hline
50 & 12.76 & 16 & 25 & 39.98 & 16 & 25 & 19.94 & 31 \\
100 & 30.29 & 15 & 50 & 51.06 & 16 & 50 & 161.66 & 32 \\ \cline{7-9}
150 & 41.91 & 16 & 75 & 15.11 & 16 & \multicolumn{3}{c}{} \\
200 & 266.45 & 16 & 100 & 259.57 & 15 & \multicolumn{3}{c}{}  \\ \cline{1-6}
\end{tabular}
\end{center}
\caption{Average performance on bilevel instances with 50 first stage variables.
Each column contains results from all 63 instances.
Each average is taken over instances where the other parameters are varied.
\label{table:averages50}}
\end{table}

The parameter choices given
in~\S\ref{sec: Cutting Planes Generation and Selection}
result in $0.3m$ disjunctive cuts, approximately $3m$ simple cuts,
and 15 bound cuts for each instance.
Also as in \S\ref{sec: Computational Result of branch-and-cut},
we experimented with not adding cutting planes in the preprocessor.
Both codes performed similarly to their respective performance with the preprocessor.
Thus, based on the results in this section
and~\S\ref{sec: Computational Result of branch-and-cut},
our default implementation is to generate cutting planes in the preprocessor.

\subsection{Inverse Quadratic Programs   \label{subsec:inverseQP}}

Jara-Moroni et al.~\cite{jaramoroni2017} presented a DC method for finding
local optima for LPCCs arising from inverse quadratic programs~\cite{jinghu3}.
The problem of interest has the form
\begin{equation}   \label{eqn.invQP}
\begin{array}{ll}
\min_{x, b, c} & || (x,b,c) - (\bar{x},\bar{b},\bar{c})||_1  \\ [5pt]
\mbox{s.t.} & x \, \in \, \mbox{argmin}_y \{ \frac{1}{2} y^TQy + c^Ty \, : \, Ay \geq b \}  \\  [5pt]
& (x,b,c) \, \in \, P
\end{array}
\end{equation}
where $\bar{x}$, $\bar{b}$, and $\bar{c}$ are observations of the parameters
and solution of a quadratic program
and $P$ is a polyhedron.
The objective is to find $(x,b,c)$ close to the observed values where
$x$ does solve the lower level quadratic program.
In our computational testing, we varied the number of rows $\tm$ and columns $\tn$
of $A$ between 100 and 400 and between 5 and 90, respectively;
the dimensions of all other vectors and matrices are determined by
the dimensions of~$A$.
When the matrix $Q$ is positive definite, the inverse QP is equivalent to
the following LPCC:
\begin{equation}   \label{eqn.invQPlpcc}
\begin{array}{ll}
\min_{x, b, c, z^x, z^b, z^c, \lambda} & \one^Tz^x + \one^Tz^b + \one^Tz^c \\
\mbox{s.t.} & Qx + c -A^T \lambda = 0  \\
& x+z^x \geq \bar{x}, \, -x+z^x \geq -\bar{x}  \\
& b+z^b \geq \bar{b}, \, -b+z^b \geq -\bar{b}  \\
& c+z^c \geq \bar{c}, \, -c+z^c \geq -\bar{c}  \\
& (x,b,c) \, \in \, P  \\
& 0 \leq \lambda \perp w := Ax -b \geq 0
\end{array}
\end{equation}
where $\lambda$ is the vector of KKT variables for the inner QP,
the variables $z^x$, $z_b$, $z^c$ are used to represent the $L_1$ objective
function in~(\ref{eqn.invQP}),
and $\mathbf{1}$ represents a vector of ones of an appropriate dimension.

The instances in \cite{jaramoroni2017} were generated in MATLAB,
whereas our instances were generated using AMPL.
Nonetheless, we closely followed their procedures except for the generation of~$Q$.
Our matrix $Q \in \real^{\tn \times \tn}$ was formed as the product $MM^T$,
where $M \in \real^{\tn \times \tn}$ was a square matrix with exactly three nonzeroes per row,
with diagonal entries uniformly distributed between 0.5 and 1
and two off-diagonal entries uniformly distributed between 0 and~1;
this results in a positive definite matrix $Q$, with about 9 entries per row
on average (similar to the number of nonzeroes in a row of $Q$ from~\cite{jaramoroni2017}).
Other parameters were generated as in~\cite{jaramoroni2017}:
the matrix $A \in \real^{\tm \times \tn}$ has an average of approximately 10 nonzero entries per row
which are uniformly distributed between 0 and 1;
a vector $\tilde{x} \in \real^{\tn}$ has components distributed as Normal(0,1);
vectors $\hat{\lambda} \in \real^{\tm}$ and $\hat{w} \in \real^{\tm}$ have components
uniformly distributed between 0 and 10;
a binary vector $v \in \bin^{\tm}$ is generated and
$\tilde{\lambda} \in \real^{\tm}$ and $\tilde{w} \in \real^{\tm}$
are constructed as the Hadamard products $\tilde{\lambda} := \hat{\lambda} \bullet v$
and $\tilde{w} := \hat{w} \bullet (\one - v)$;
vectors $\tilde{b} \in \real^{\tm}$ and $\tilde{c} \in \real^{\tn}$ are defined as
$\tilde{b}:=A\tilde{x}-\tilde{w}$ and $\tilde{c}=A^T\tilde{\lambda}-Q\tilde{x}$;
vectors $\bar{x} \in \real^{\tn}$, $\bar{b} \in \real^{\tm}$, and $\bar{c} \in \real^{\tn}$ are obtained by
perturbing $\tilde{x}$, $\tilde{b}$, and $\tilde{c}$ respectively,
using Normal (0,1) noise;
the polyhedron $P$ is constructed as a box using
simple bounds  $-u^x \one \leq x \leq u^x \one$,
$-u^b \one \leq b \leq u^b \one$, $-u^c \one \leq c \leq u^c \one$ with
$u^x=10 \max \{|\tilde{x}_i|\}$,
$u^b=10 \max \{|\tilde{b}_i|\}$,
$u^c=10 \max \{|\tilde{c}_i|\}$;
finally, upper bounds are also imposed on $\lambda$ with
$u^{\lambda}=10 \max \{|\tilde{\lambda}_i|\}$.
The point $(\tilde{x},\tilde{b},\tilde{c})$
with $\tilde{\lambda}$ is feasible in the resulting problem instances
of~(\ref{eqn.invQPlpcc}).

It is easy to generate explicit upper bounds on $w=Ax-b$
from the upper bounds on $x$ and~$b$.
Also, explicit upper bounds on $\lambda$ are imposed following~\cite{jaramoroni2017}.
Thus, this problem can be formulated directly
as a mixed integer program of the form (\ref{eq:MIP LPCC}).
Because of this observation, our comparisons in this section
are somewhat different from the previous experiments.
In particular, we make the following two changes:
\begin{itemize}
\item Since bounds are already available, we do not use the  cutting plane
generation features of the preprocessor.
\item We compare our LPCC branch-and-cut code with the CPLEX MIP solver
invoked from AMPL, run with a single thread.
\end{itemize}

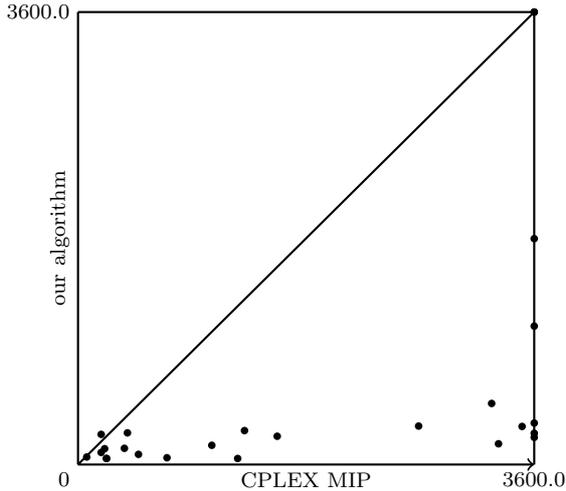
\begin{figure}
\begin{center}
\begin{tikzpicture}[scale=0.5]

\draw[thick] (0,12) -- (12,12) -- (12,0);
\draw[thick] (0,0) -- (12,12);
\draw[thick] (0,0) -- (0,12) 
  node[midway, sloped,above] {our algorithm};
\draw[thick,->] (0,0) -- (12,0) node[midway, below] {CPLEX MIP};

\draw (0,0) node [below left] { 0 };
\draw (12,0) node [below] { 3600.0 };

\draw (0,12) node [left] { 3600.0 };

\fill (11.68,1.01) circle(0.1);
\fill (8.96,1.02) circle(0.1);
\fill (0.23,0.2) circle(0.1);
\fill (12.0,1.1) circle(0.1);
\fill (0.61,0.32) circle(0.1);
\fill (4.38,0.9) circle(0.1);
\fill (12.0,0.83) circle(0.1);
\fill (0.7,0.42) circle(0.1);
\fill (12.0,0.72) circle(0.1);
\fill (12.0,3.67) circle(0.1);
\fill (2.34,0.18) circle(0.1);
\fill (10.88,1.62) circle(0.1);
\fill (11.06,0.55) circle(0.1);
\fill (4.2,0.16) circle(0.1);
\fill (12.0,12.0) circle(0.1);
\fill (3.52,0.51) circle(0.1);
\fill (1.59,0.27) circle(0.1);
\fill (12.0,12.0) circle(0.1);
\fill (1.22,0.43) circle(0.1);
\fill (0.75,0.16) circle(0.1);
\fill (5.24,0.75) circle(0.1);
\fill (12.0,5.99) circle(0.1);
\fill (0.61,0.8) circle(0.1);
\fill (1.3,0.84) circle(0.1);
\fill (0.75,0.16) circle(0.1);

\end{tikzpicture}
\end{center}
\caption{\label{fig:scatter_invQP}Scatter plots for CPU time (seconds) for solution of inverse quadratic
programs.
Horizontal axis is time for CPLEX MIP called from AMPL and run on a single thread,
vertical axis is time for our branch and cut algorithm.
}
\end{figure}

Our testbed consisted of 5 sets of 5 instances:
$(\tm,\tn)$ equal to (100,75), (120,90), (150,20), (200,15), and (400,5). 
A scatter plot of the results can be found in Figure~\ref{fig:scatter_invQP}
and performance profiles can be found in Figure~\ref{fig:perfprof_invQP}.
Detailed computational results are contained in the Appendix,
in Table~\ref{table.inverseQP}.
Our algorithm was able to solve 23 of the 25 instances within
the 3600 second time limit;
the corresponding figure for CPLEX was 18 out of 25.
There was only one instance where CPLEX outperformed our code.
Our algorithm solved 20 of the 25 instances within 360 seconds,
while CPLEX only solved 6 of the instances within this time window.

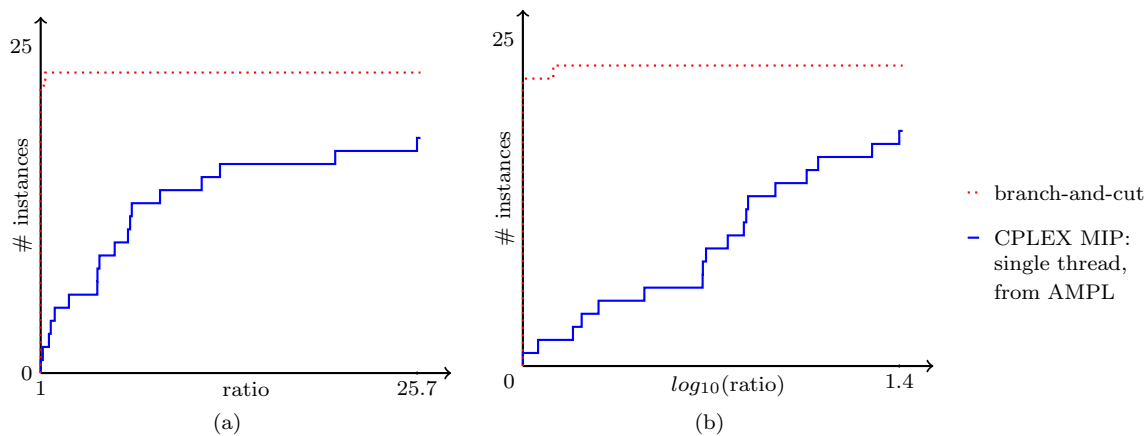
\begin{figure}
\begin{center}
\begin{tabular}{ccl}
\begin{tikzpicture}[scale=0.45]

\draw (11,0) node [below] { 25.7 };

\draw (11,-0.05) -- (11,0.05);
\draw (0,0) node [below] { 1 };
\draw (0,0) node [left] { 0 };
\draw (0,9.625) node [left]  { 25 };

\draw[thick,->] (0,0) -- (0,10.5)
  node[midway, sloped,above] {$\#$ instances};

\draw[thick,->] (0,0) -- (12,0) node[midway, below] {ratio};
\draw[thick, blue] (0,0) -- (0.0,0.385) -- (0.06243999999999996,0.385) -- (0.06243999999999996,0.77) -- (0.24084000000000003,0.77) -- (0.24084000000000003,1.155) -- (0.29435999999999996,1.155) -- (0.29435999999999996,1.54) -- (0.41031999999999996,1.54) -- (0.41031999999999996,1.925) -- (0.8251000000000001,1.925) -- (0.8251000000000001,2.31) -- (1.65466,2.31) -- (1.65466,2.6950000000000003) -- (1.65912,2.6950000000000003) -- (1.65912,3.08) -- (1.7170999999999998,3.08) -- (1.7170999999999998,3.465) -- (2.1631,3.465) -- (2.1631,3.85) -- (2.55112,3.85) -- (2.55112,4.235) -- (2.61356,4.235) -- (2.61356,4.62) -- (2.66262,4.62) -- (2.66262,5.005) -- (3.4877200000000004,5.005) -- (3.4877200000000004,5.390000000000001) -- (4.7053,5.390000000000001) -- (4.7053,5.775) -- (5.2405,5.775) -- (5.2405,6.16) -- (8.607800000000001,6.16) -- (8.607800000000001,6.545) -- (10.99836,6.545) -- (10.99836,6.93) -- (11.1,6.93);

\draw[thick, red, dotted] (0,0) -- (0.0,0.385) -- (0.0,0.77) -- (0.0,1.155) -- (0.0,1.54) -- (0.0,1.925) -- (0.0,2.31) -- (0.0,2.6950000000000003) -- (0.0,3.08) -- (0.0,3.465) -- (0.0,3.85) -- (0.0,4.235) -- (0.0,4.62) -- (0.0,5.005) -- (0.0,5.390000000000001) -- (0.0,5.775) -- (0.0,6.16) -- (0.0,6.545) -- (0.0,6.93) -- (0.0,7.315) -- (0.0,7.7) -- (0.0,8.085) -- (0.0,8.47) -- (0.13380000000000003,8.47) -- (0.13380000000000003,8.855) -- (11.1,8.855);

\end{tikzpicture}

&

\begin{tikzpicture}[scale=0.45]

\draw (11,0) node [below] { 1.4 };

\draw (11,-0.05) -- (11,0.05);
\draw (0,0) node [below left] { 0 };
\draw (0,9.625) node [left]  { 25 };

\draw[thick,->] (0,0) -- (0,10.5)
  node[midway, sloped,above] {$\#$ instances};

\draw[thick,->] (0,0) -- (12,0) node[midway, below] {$log_{10}$(ratio)};
\draw[thick, blue] (0,0) -- (0.0,0.385) -- (0.44419926953250477,0.385) -- (0.44419926953250477,0.77) -- (1.4637867468494308,0.77) -- (1.4637867468494308,1.155) -- (1.71816373524067,1.155) -- (1.71816373524067,1.54) -- (2.211449391259908,1.54) -- (2.211449391259908,1.925) -- (3.5505189772264307,1.925) -- (3.5505189772264307,2.31) -- (5.253601201048164,2.31) -- (5.253601201048164,2.6950000000000003) -- (5.260791241337689,2.6950000000000003) -- (5.260791241337689,3.08) -- (5.35290001152927,3.08) -- (5.35290001152927,3.465) -- (5.9884186906375,3.465) -- (5.9884186906375,3.85) -- (6.45844454545816,3.85) -- (6.45844454545816,4.235) -- (6.528346047206904,4.235) -- (6.528346047206904,4.62) -- (6.582275065833062,4.62) -- (6.582275065833062,5.005) -- (7.380327775538985,5.005) -- (7.380327775538985,5.390000000000001) -- (8.294514968885041,5.390000000000001) -- (8.294514968885041,5.775) -- (8.629612502314417,5.775) -- (8.629612502314417,6.16) -- (10.20631407195054,6.16) -- (10.20631407195054,6.545) -- (11.000657425815728,6.545) -- (11.000657425815728,6.93) -- (11.1,6.93);

\draw[thick, red, dotted] (0,0) -- (0.0,0.385) -- (0.0,0.77) -- (0.0,1.155) -- (0.0,1.54) -- (0.0,1.925) -- (0.0,2.31) -- (0.0,2.6950000000000003) -- (0.0,3.08) -- (0.0,3.465) -- (0.0,3.85) -- (0.0,4.235) -- (0.0,4.62) -- (0.0,5.005) -- (0.0,5.390000000000001) -- (0.0,5.775) -- (0.0,6.16) -- (0.0,6.545) -- (0.0,6.93) -- (0.0,7.315) -- (0.0,7.7) -- (0.0,8.085) -- (0.0,8.47) -- (0.889441808107168,8.47) -- (0.889441808107168,8.855) -- (11.1,8.855);

\end{tikzpicture}
&

\begin{tikzpicture}[scale=0.3]
\draw[thick, red, dotted] (0,9) -- (0.5,9);
\draw (0.8,9) node [right] {branch-and-cut};
\draw[thick, blue] (0,7) -- (0.5,7);
\draw (0.8,7) node [right] {CPLEX MIP:};
\draw (0.8,5.8) node [right] {single thread,};
\draw (0.8,4.6) node [right] {from AMPL};
\draw (0,0) node [right] {$\quad$};
\end{tikzpicture}

\\
(a) & (b)
\end{tabular}
\end{center}
\caption{\label{fig:perfprof_invQP}Performance profile for CPU time (seconds) for solution of
25 inverse quadratic programs.
Vertical axis is the number of instances.
Horizontal axis is ratio of time required by the given algorithm to the
time required by the better algorithm.
Time limit 3600 seconds.
(a) Linear scale.  (b) Log scale.}
\end{figure}

\section{Conclusions}

The optimal solution to a linear program with complementarity constraints
can in principle be found directly using CPLEX.
However, far better performance can often be obtained
by adding good cutting planes, by incorporating a specialized feasibility
recovery routine, and especially by designing good branching routines.
Our computational results demonstrate that our code is at least
an order of magnitude faster than a default version of CPLEX,
at least for our test set of instances.
It is able to solve instances with up to 400 complementarity constraints in reasonable amounts
of time,
and can reliably solve instances with 100 complementarity constraints
in less than a minute.

\bibliography{optim}

\appendix

\section{LPCC Test Instances  \label{appendix:test instances}}
In order to test the effectiveness of different type of valid constraints, a series of LPCC instances was randomly generated, and Procedure \ref{proc: LPCC instances generator 1} gives a detailed description of the generator. 

\begin{algorithm}
\SetAlgorithmName{Procedure}{}

\SetKwInOut{Input}{input}\SetKwInOut{Output}{output}

\Input{$n$,$m$,$k$,$rankM$,$dense$}
\Output{vector $c$,$d$,$b$,$q$; matrix $A$,$B$,$N$,$M$ }
\BlankLine
\textit{1:} generate $n$ dimension vector $\bar{x}$ with value between 0 and 10, integer\;

\textit{2:} generate $m$ dimension vector $\bar{y}$ with value between 0 and 10, integer if index $<\dfrac{m}{3}$; 0 otherwise\;
\textit{3:} generate $n$ dimension vector $c$ with value between 0 and 10, integer\;
\textit{4:} generate $m$ dimension  vector $d$ with value between 0 and 10, integer\;
\textit{5:} generate $k\times n$ matrix $A$ with value between -5 and 6, integer, and the matrix density is $dense$\;
\textit{6:} generate $k\times m$ matrix $B$ with value between -5 and 6, integer, and the matrix density is $dense$\;
\textit{7:} generate $m\times n$ matrix $N$ with value between -5 and 6, integer, and the matrix density is $dense$\;
\textit{8:} generate $m\times rankM$ matrix $L$ with value between -5 and 6, integer, and the matrix density is $dense$; generate $m \times m$ upper triangular matrix $\Delta M$ with value between -2 and 2, integer; Let $m \times m$ matrix $M=LL^T+\Delta M-\Delta M^T$\;
\textit{9:} generate $k$ dimension vector $\Delta b$ with value between 1 and 11, integer; let $k$ dimension vector $b=A\bar{x}+B\bar{y}-\Delta b$\;
\textit{10:} generate $m$ dimension vector $\Delta q$ with value 0 if index $<\dfrac{2m}{3}$; integer between 1 and 11 otherwise; let $m$ dimension vector $q=-N\bar{x}-M\bar{y}+\Delta q$\;
\caption{LPCC instances generator}\label{proc: LPCC instances generator 1}
\end{algorithm}
\begin{remark}
In the initialization step of the procedure,
$n$ is the dimension of $x$ variable; $m$ is the dimension of $y$ variable; $k$ is the dimension of $b$;
$rankM$ is the rank of matrix $M$; $dense$ is the density of generated matrices;
we assume all instances have the non-negativity constraint $x\geq0$ which are not included in the constraint $Ax+By\geq b$; \textit{step 1} and \textit{step 2} are used to generate a feasible LPCC solution; step 8 is to generate matrix $M$ to be a non-symmetric positive semidefinite matrix with rank $rankM$. 
\end{remark}

We  generated 60 LPCC instances with 100, 150, 200 complementaries, 20 instances of each size, and with the same parameter, we randomly generated 5 instances.
For CPLEX solving LPCC instances, we used indicator constraints in CPLEX C callable library \cite{ilog5} to formulate the complementarity constraints, and the CPLEX setting is default. The time limit for CPLEX is 3600 seconds. Notice that default CPLEX is unable to solve most of our LPCC instances when $m=200$ within 3600 seconds.
Table \ref{table:data60LPCC} contains objective function value information for the 60 instances,
including the effectiveness of the preprocessing routines.

\begin{table}
\begin{center}
{
\begin{tabular}{|r|r|r|r|r|r|r|r|r|r|r|}
\hline
&&&& Optimal & LP & \multicolumn{2}{c|}{Preprocessed bounds} & \multicolumn{3}{c|}{Relative gaps (percentages)}
\\  \cline{7-11}
\#	&	m	&	rank	&	dense	&	Value	&	relaxation	&	lower	&	upper	&	UB-LB	&	UB-opt	&	\% closed		\\ \hline
1	&	100	&	30	&	70	&	769.911528	&	629.002874	&	669.22439	&	770.287	&	13.13	&	0.05	&	71.46		\\
2	&	100	&	30	&	70	&	752	&	650.929154	&	723.79249	&	754.658	&	4.10	&	0.35	&	27.91		\\
3	&	100	&	30	&	70	&	690.306012	&	627.332027	&	657.571025	&	691	&	4.84	&	0.10	&	51.98		\\
4	&	100	&	30	&	70	&	543	&	531.188245	&	539.856029	&	544.497	&	0.85	&	0.28	&	26.62		\\
5	&	100	&	30	&	70	&	930	&	771.820799	&	896.57115	&	930.917	&	3.69	&	0.10	&	21.13		\\
6	&	100	&	30	&	20	&	589	&	583.487434	&	588.868287	&	589	&	0.02	&	0.00	&	2.39		\\
7	&	100	&	30	&	20	&	488	&	425.717966	&	459.942655	&	488	&	5.75	&	0.00	&	45.05		\\
8	&	100	&	30	&	20	&	771	&	687.744893	&	745.909078	&	771	&	3.25	&	0.00	&	30.14		\\
9	&	100	&	30	&	20	&	628	&	524.270776	&	620.866694	&	628	&	1.14	&	0.00	&	6.88		\\
10	&	100	&	30	&	20	&	732	&	705.051229	&	729.547378	&	732	&	0.34	&	0.00	&	9.10		\\
11	&	100	&	60	&	70	&	612.145738	&	606.45432	&	609.833638	&	622.283	&	2.03	&	1.66	&	40.62		\\
12	&	100	&	60	&	70	&	686.130259	&	649.068458	&	675.208522	&	686.212	&	1.60	&	0.01	&	29.47		\\
13	&	100	&	60	&	70	&	734	&	722.033536	&	733.174515	&	734	&	0.11	&	0.00	&	6.90		\\
14	&	100	&	60	&	70	&	665.868588	&	657.703283	&	661.460391	&	666	&	0.68	&	0.02	&	53.99		\\
15	&	100	&	60	&	70	&	984.588193	&	818.248599	&	855.932906	&	986	&	13.21	&	0.14	&	77.34		\\
16	&	100	&	60	&	20	&	691	&	629.620621	&	664.054558	&	691	&	3.90	&	0.00	&	43.90		\\
17	&	100	&	60	&	20	&	666.995818	&	631.110603	&	655.171515	&	667	&	1.77	&	0.00	&	32.95		\\
18	&	100	&	60	&	20	&	756.780603	&	725.103749	&	746.527684	&	758	&	1.52	&	0.16	&	32.37		\\
19	&	100	&	60	&	20	&	763	&	626.529227	&	722.010148	&	763.971	&	5.50	&	0.13	&	30.04		\\
20	&	100	&	60	&	20	&	532.218697	&	521.894551	&	528.196096	&	533	&	0.90	&	0.15	&	38.96		\\  \hline
21	&	150	&	30	&	70	&	1029	&	946.929565	&	1010.002422	&	1029	&	1.85	&	0.00	&	23.15		\\
22	&	150	&	30	&	70	&	1160	&	1075.719667	&	1143.215912	&	1160	&	1.45	&	0.00	&	19.91		\\
23	&	150	&	30	&	70	&	965	&	929.722695	&	957.060812	&	965	&	0.82	&	0.00	&	22.51		\\
24	&	150	&	30	&	70	&	1242	&	1170.744571	&	1232.634488	&	1242	&	0.75	&	0.00	&	13.14		\\
25	&	150	&	30	&	70	&	1149	&	1013.045865	&	1063.947928	&	1149	&	7.40	&	0.00	&	62.56		\\
26	&	150	&	30	&	20	&	822.333333	&	790.161133	&	813.932095	&	822.333	&	1.02	&	0.00	&	26.11		\\
27	&	150	&	30	&	20	&	1046	&	991.351886	&	1039.478766	&	1046	&	0.62	&	0.00	&	11.93		\\
28	&	150	&	30	&	20	&	922	&	851.085225	&	899.489258	&	922	&	2.44	&	0.00	&	31.74		\\
29	&	150	&	30	&	20	&	992	&	855.028214	&	921.051941	&	992	&	7.15	&	0.00	&	51.80		\\
30	&	150	&	30	&	20	&	848	&	729.617101	&	775.254605	&	848	&	8.58	&	0.00	&	61.45		\\
31	&	150	&	100	&	70	&	1377.072388	&	1263.798462	&	1344.135656	&	1377.072	&	2.39	&	0.00	&	29.08		\\
32	&	150	&	100	&	70	&	837	&	833.238632	&	835.993215	&	837	&	0.12	&	0.00	&	26.77		\\
33	&	150	&	100	&	70	&	972.779519	&	912.297933	&	951.089989	&	972.804	&	2.23	&	0.00	&	35.86		\\
34	&	150	&	100	&	70	&	1260.57242	&	1206.833191	&	1238.300018	&	1261.188	&	1.82	&	0.05	&	41.45		\\
35	&	150	&	100	&	70	&	1087.08492	&	1040.170448	&	1077.111477	&	1089	&	1.09	&	0.18	&	21.26		\\
36	&	150	&	100	&	20	&	921.273479	&	893.518557	&	904.290053	&	923	&	2.03	&	0.19	&	61.19		\\
37	&	150	&	100	&	20	&	923.772654	&	774.71571	&	879.664636	&	925	&	4.91	&	0.13	&	29.59		\\
38	&	150	&	100	&	20	&	1139	&	1111.79451	&	1126.884941	&	1139	&	1.06	&	0.00	&	44.53		\\
39	&	150	&	100	&	20	&	879.582356	&	812.660589	&	852.096526	&	879.605	&	3.13	&	0.00	&	41.07		\\
40	&	150	&	100	&	20	&	1158.383138	&	1063.017814	&	1119.548217	&	1158.432	&	3.36	&	0.00	&	40.72		\\  \hline
41	&	200	&	30	&	70	&	1580	&	1098.044624	&	1196.5995	&	1580	&	24.27	&	0.00	&	79.55		\\
42	&	200	&	30	&	70	&	1057	&	1025.39776	&	1050.433637	&	1057	&	0.62	&	0.00	&	20.78		\\
43	&	200	&	30	&	70	&	1577	&	1467.609941	&	1541.862973	&	1577	&	2.23	&	0.00	&	32.12		\\
44	&	200	&	30	&	70	&	1535	&	1462.36974	&	1524.019988	&	1535	&	0.72	&	0.00	&	15.12		\\
45	&	200	&	30	&	70	&	1153	&	1122.856763	&	1145.503839	&	1153	&	0.65	&	0.00	&	24.87		\\
46	&	200	&	30	&	20	&	1229	&	1148.301545	&	1192.605532	&	1229	&	2.96	&	0.00	&	45.10		\\
47	&	200	&	30	&	20	&	1350	&	1251.324462	&	1318.973919	&	1350	&	2.30	&	0.00	&	31.44		\\
48	&	200	&	30	&	20	&	1451	&	1115.387691	&	1208.887517	&	1451	&	16.69	&	0.00	&	72.14		\\
49	&	200	&	30	&	20	&	1345	&	1261.123305	&	1337.276135	&	1345	&	0.57	&	0.00	&	9.21		\\
50	&	200	&	30	&	20	&	1249	&	1164.340236	&	1195.763472	&	1249	&	4.26	&	0.00	&	62.88		\\
51	&	200	&	120	&	70	&	1726.526853	&	1649.267937	&	1701.696186	&	1728	&	1.52	&	0.09	&	32.14		\\
52	&	200	&	120	&	70	&	1403	&	1337.168109	&	1394.142467	&	1403	&	0.63	&	0.00	&	13.45		\\
53	&	200	&	120	&	70	&	1144.989488	&	1126.310832	&	1143.367197	&	1145	&	0.14	&	0.00	&	8.69		\\
54	&	200	&	120	&	70	&	1542	&	1500.576683	&	1532.123758	&	1542	&	0.64	&	0.00	&	23.84		\\
55	&	200	&	120	&	70	&	1096.255705	&	951.763018	&	1009.459289	&	1097	&	7.99	&	0.07	&	60.07		\\
56	&	200	&	120	&	20	&	1235.593203	&	1183.04243	&	1203.799741	&	1237	&	2.69	&	0.11	&	60.50		\\
57	&	200	&	120	&	20	&	1224.764683	&	1100.94521	&	1188.734874	&	1226	&	3.04	&	0.10	&	29.10		\\
58	&	200	&	120	&	20	&	1145.969792	&	1132.996319	&	1140.093917	&	1147	&	0.60	&	0.09	&	45.29		\\
59	&	200	&	120	&	20	&	1426	&	1399.225251	&	1415.85159	&	1429.364	&	0.95	&	0.24	&	37.90		\\
60	&	200	&	120	&	20	&	1371.901959	&	1340.784415	&	1358.035244	&	1372	&	1.02	&	0.01	&	44.56		\\ \hline
\multicolumn{8}{r|}{Means:} & 3.28 & 0.07 & 35.40 \\ \cline{9-11}
\end{tabular}
}
\end{center}
\caption{Objective function data for the 60 instances.  \label{table:data60LPCC}
All instances have $n=2$ and $k=20$.
The number of complementarities is $m$.
The rank of $M$ and the density of each matrix are indicated.
Three relative gaps are given as percentages:
(i) the gap between the upper and lower bounds obtained through preprocessing,
(ii) the gap between the upper bound obtained from feasibility recovery and
the optimal value of the LPCC,
and (iii) the improvement in the gap between upper and lower bound effected by
the improvement in the LP relaxation obtained through preprocessing.}
\end{table}
Table \ref{table:results60LPCC} contains
performance data.
\begin{table}
\begin{center}
{
\begin{tabular}{|r|r|r|r|r|r|r|r|r|r|}
\hline
&&&& Preprocess & \multicolumn{2}{c|}{Our algorithm} & \multicolumn{3}{c|}{default CPLEX indicator constraint} \\ \cline{6-10}
\#	&	m	&	rank	&	dense	&		time	&	time	&	nodes & time 	& nodes 	&	\% gap		\\ \hline
1	&	100	&	30	&	70	&		5.93	&	1.19	&	404	&		7.04	&	2704	&			\\
2	&	100	&	30	&	70	&		4.71	&	5.70	&	4340	&		DNF	&	6572052	&	0.485		\\
3	&	100	&	30	&	70	&		9.35	&	0.61	&	260	&		4.98	&	1775	&			\\
4	&	100	&	30	&	70	&		1.79	&	0.13	&	38	&		1.70	&	185	&			\\
5	&	100	&	30	&	70	&		4.80	&	0.46	&	210	&		4.18	&	1330	&			\\
6	&	100	&	30	&	20	&		1.02	&	0.02	&	18	&		0.73	&	43	&			\\
7	&	100	&	30	&	20	&		4.52	&	0.24	&	44	&		0.81	&	103	&			\\
8	&	100	&	30	&	20	&		5.55	&	0.55	&	224	&		5.24	&	3305	&			\\
9	&	100	&	30	&	20	&		3.56	&	0.08	&	32	&		1.04	&	204	&			\\
10	&	100	&	30	&	20	&		0.81	&	0.10	&	30	&		1.16	&	208	&			\\
11	&	100	&	60	&	70	&		1.30	&	0.24	&	88	&		7.08	&	3940	&			\\
12	&	100	&	60	&	70	&		6.38	&	0.73	&	588	&		384.44	&	353435	&			\\
13	&	100	&	60	&	70	&		1.21	&	0.10	&	20	&		1.09	&	49	&			\\
14	&	100	&	60	&	70	&		1.21	&	1.06	&	534	&		9.32	&	4671	&			\\
15	&	100	&	60	&	70	&		7.25	&	3.69	&	1322	&		DNF	&	9693517	&	0.125		\\
16	&	100	&	60	&	20	&		5.29	&	1.01	&	444	&		2.91	&	907	&			\\
17	&	100	&	60	&	20	&		6.02	&	2.06	&	1324	&		7.47	&	5723	&			\\
18	&	100	&	60	&	20	&		1.21	&	0.35	&	206	&		5.55	&	2672	&			\\
19	&	100	&	60	&	20	&		4.76	&	0.37	&	194	&		3.73	&	1011	&			\\
20	&	100	&	60	&	20	&		1.00	&	0.37	&	226	&		2.83	&	708	&			\\  \hline
\multicolumn{4}{r|}{\mbox{Means:}} & 3.88 & 0.95 & 527 \\ \hline
21	&	150	&	30	&	70	&		24.66	&	2.40	&	448	&		15.58	&	3307	&			\\
22	&	150	&	30	&	70	&		23.08	&	0.96	&	124	&		4.27	&	0	&			\\
23	&	150	&	30	&	70	&		5.11	&	0.16	&	56	&		4.58	&	211	&			\\
24	&	150	&	30	&	70	&		24.49	&	0.61	&	98	&		8.52	&	2163	&			\\
25	&	150	&	30	&	70	&		24.69	&	6.17	&	942	&		685.24	&	207429	&			\\
26	&	150	&	30	&	20	&		3.68	&	0.20	&	92	&		3.27	&	232	&			\\
27	&	150	&	30	&	20	&		16.81	&	0.20	&	32	&		2.76	&	124	&			\\
28	&	150	&	30	&	20	&		18.82	&	0.36	&	78	&		3.62	&	370	&			\\
29	&	150	&	30	&	20	&		14.23	&	1.60	&	188	&		7.03	&	1256	&			\\
30	&	150	&	30	&	20	&		18.49	&	3.70	&	682	&		22.30	&	10877	&			\\
31	&	150	&	100	&	70	&		26.75	&	101.60	&	28538	&		DNF	&	13175000	&	0.146		\\
32	&	150	&	100	&	70	&		4.86	&	0.54	&	192	&		6.31	&	580	&			\\
33	&	150	&	100	&	70	&		27.11	&	29.50	&	8744	&		DNF	&	2323005	&	0.114		\\
34	&	150	&	100	&	70	&		14.16	&	132.03	&	32124	&		DNF	&	1674151	&	0.206		\\
35	&	150	&	100	&	70	&		5.21	&	10.90	&	2888	&		DNF	&	2443005	&	0.272		\\
36	&	150	&	100	&	20	&		4.42	&	16.12	&	4602	&		DNF	&	15955452	&	0.437		\\
37	&	150	&	100	&	20	&		22.16	&	15.45	&	3064	&		2338.14	&	845218	&			\\
38	&	150	&	100	&	20	&		4.84	&	2.32	&	584	&		11.66	&	1427	&			\\
39	&	150	&	100	&	20	&		22.38	&	4.84	&	976	&		19.74	&	4403	&			\\
40	&	150	&	100	&	20	&		23.75	&	32.48	&	10376	&		2063.57	&	1202357	&			\\  \hline
\multicolumn{4}{r|}{\mbox{Means:}} & 15.49 & 16.53 & 4249 \\ \hline
41	&	200	&	30	&	70	&		126.63	&	1546.32	&	181008	&		DNF	&	1368269	&	0.612		\\
42	&	200	&	30	&	70	&		12.97	&	0.22	&	38	&		5.50	&	0	&			\\
43	&	200	&	30	&	70	&		76.18	&	1.19	&	66	&		8.43	&	189	&			\\
44	&	200	&	30	&	70	&		15.66	&	1.48	&	262	&		29.50	&	3328	&			\\
45	&	200	&	30	&	70	&		14.40	&	0.33	&	64	&		5.77	&	0	&			\\
46	&	200	&	30	&	20	&		44.60	&	0.85	&	42	&		5.51	&	0	&			\\
47	&	200	&	30	&	20	&		49.62	&	1.17	&	120	&		15.08	&	1714	&			\\
48	&	200	&	30	&	20	&		58.34	&	22.07	&	930	&		1538.63	&	405689	&			\\
49	&	200	&	30	&	20	&		39.58	&	0.69	&	48	&		5.75	&	0	&			\\
50	&	200	&	30	&	20	&		48.30	&	2.80	&	196	&		19.23	&	2749	&			\\
51	&	200	&	120	&	70	&		13.90	&	176.16	&	27046	&		DNF	&	748850	&	0.292		\\
52	&	200	&	120	&	70	&		14.87	&	319.53	&	53786	&		DNF	&	635814	&	0.048		\\
53	&	200	&	120	&	70	&		15.17	&	25.74	&	5014	&		DNF	&	1721873	&	0.036		\\
54	&	200	&	120	&	70	&		13.97	&	47.07	&	7736	&		DNF	&	1350106	&	0.028		\\
55	&	200	&	120	&	70	&		87.88	&	86.08	&	8646	&		DNF	&	4416595	&	0.168		\\
56	&	200	&	120	&	20	&		14.61	&	283.36	&	42010	&		DNF	&	1371796	&	0.248		\\
57	&	200	&	120	&	20	&		82.14	&	1017.50	&	97688	&		DNF	&	651199	&	0.736		\\
58	&	200	&	120	&	20	&		12.90	&	11.19	&	1834	&		DNF	&	1233196	&	0.160		\\
59	&	200	&	120	&	20	&		14.43	&	31.69	&	5188	&		DNF	&	856913	&	0.101		\\
60	&	200	&	120	&	20	&		13.87	&	491.75	&	85978	&		DNF	&	710977	&	0.278		\\ \hline
\multicolumn{4}{r|}{\mbox{Means:}} & 58.50 & 203.36 & 25885 \\ \cline{5-7}
\end{tabular}
}
\end{center}
\caption{Performance data for the 60 instances.   \label{table:results60LPCC}
The final gap obtained by default CPLEX is indicated for the 18 instances it didn't solve (DNF)
within the time limit of 3600 seconds.}
\end{table}

\section{Bilevel Test Instances   \label{app:bilevel}}

Our code solved all 63 of the instances with dimension of $v$
equal to 50 and $18/35$ of the larger instances.
With extended time, default CPLEX was able to solve all but one problem with
$n=50$; it still has a gap of $16.56$\% for problem 60 after more than 7200 seconds of
wall clock time and 47304 seconds of processor time.
It solved just $6/35$ of the larger instances.
Run time information can be found in Tables \ref{table:perf_bil_50} and~\ref{table:perf_bil_larger}.
\begin{table}
\begin{center}
{\scriptsize
\begin{tabular}{|r|r|r|r|r|r|r|r|r|}
\hline
& \multicolumn{3}{c|}{Dimensions} && LP & Preprocess & Optimal & \% Gap \\ \cline{2-4}
\#	&	$v$	&	$b$	&	$g$	&	rank($Q$)	&	relaxation	&	lower bound	&	value	& shrunk	\\ \hline
1	&	50	&	25	&	50	&	25	&	0.644247	&	0.67043	&	0.708016	&	41.06					\\
2	&	50	&	25	&	50	&	25	&	0.691488	&	0.742312	&	0.817536	&	40.32					\\
3	&	50	&	25	&	100	&	25	&	0.758156	&	0.853512	&	0.90403	&	65.37					\\
4	&	50	&	25	&	100	&	25	&	0.508003	&	0.719332	&	1.030281	&	40.46					\\
5	&	50	&	25	&	150	&	25	&	0.634865	&	0.768284	&	0.930512	&	45.13					\\
6	&	50	&	25	&	150	&	25	&	0.746323	&	0.939574	&	1.179523	&	44.61					\\
7	&	50	&	25	&	200	&	25	&	0.834849	&	0.882683	&	0.983496	&	32.18					\\
8	&	50	&	25	&	200	&	25	&	0.552456	&	0.742789	&	1.113869	&	33.90					\\
\hline
9	&	50	&	50	&	50	&	25	&	0.74926	&	0.833068	&	0.974266	&	37.25					\\
10	&	50	&	50	&	50	&	25	&	0.596532	&	0.680096	&	0.777676	&	46.13					\\
11	&	50	&	50	&	100	&	25	&	0.74352	&	0.85814	&	1.025396	&	40.66					\\
12	&	50	&	50	&	100	&	25	&	0.713623	&	0.897497	&	1.007106	&	62.65					\\
13	&	50	&	50	&	150	&	25	&	0.734739	&	0.788571	&	0.884411	&	35.97					\\
14	&	50	&	50	&	150	&	25	&	0.521193	&	0.682435	&	0.990508	&	34.36					\\
15	&	50	&	50	&	200	&	25	&	0.754768	&	0.782534	&	0.783561	&	96.43					\\
16	&	50	&	50	&	200	&	25	&	0.754747	&	0.847708	&	1.054281	&	31.04					\\
\hline
17	&	50	&	75	&	50	&	25	&	0.649778	&	0.773429	&	0.945024	&	41.88					\\
18	&	50	&	75	&	50	&	25	&	0.607476	&	0.880724	&	1.046959	&	62.17					\\
19	&	50	&	75	&	100	&	25	&	0.720769	&	0.812158	&	0.971363	&	36.47					\\
20	&	50	&	75	&	100	&	25	&	0.529814	&	0.667117	&	0.949297	&	32.73					\\
21	&	50	&	75	&	150	&	25	&	0.909594	&	0.93072	&	0.933821	&	87.20					\\
22	&	50	&	75	&	150	&	25	&	0.710307	&	0.836857	&	1.103089	&	32.22					\\
23	&	50	&	75	&	200	&	25	&	0.718915	&	0.979733	&	1.326493	&	42.93					\\
24	&	50	&	75	&	200	&	25	&	0.794803	&	0.916565	&	0.950861	&	78.02					\\ 
\hline
25	&	50	&	100	&	50	&	25	&	0.766494	&	0.894661	&	1.086423	&	40.06					\\
26	&	50	&	100	&	50	&	25	&	0.485909	&	0.627154	&	0.962494	&	29.64					\\
27	&	50	&	100	&	100	&	25	&	0.767284	&	0.838957	&	0.95196	&	38.81					\\
28	&	50	&	100	&	150	&	25	&	0.578038	&	0.699594	&	0.793066	&	56.53					\\
29	&	50	&	100	&	150	&	25	&	0.713984	&	0.743964	&	0.760946	&	63.84					\\
30	&	50	&	100	&	200	&	25	&	0.616827	&	0.867557	&	1.064487	&	56.01					\\
31	&	50	&	100	&	200	&	25	&	0.651844	&	0.751472	&	0.793559	&	70.30					\\
\hline
32	&	50	&	25	&	50	&	50	&	0.644746	&	0.770222	&	0.897356	&	49.67					\\
33	&	50	&	25	&	50	&	50	&	0.602	&	0.759655	&	0.928742	&	48.25					\\
34	&	50	&	25	&	100	&	50	&	0.660691	&	0.816037	&	1.03487	&	41.52					\\
35	&	50	&	25	&	100	&	50	&	0.578159	&	0.804343	&	1.031041	&	49.94					\\
36	&	50	&	25	&	150	&	50	&	0.69423	&	0.856205	&	1.065649	&	43.61					\\
37	&	50	&	25	&	150	&	50	&	0.790053	&	0.903077	&	1.01511	&	50.22					\\
38	&	50	&	25	&	200	&	50	&	0.619887	&	0.766884	&	0.977958	&	41.05					\\
39	&	50	&	25	&	200	&	50	&	0.661197	&	0.876876	&	1.104905	&	48.61					\\
\hline
40	&	50	&	50	&	50	&	50	&	0.57275	&	0.793187	&	1.07134	&	44.21					\\
41	&	50	&	50	&	50	&	50	&	0.54549	&	0.696252	&	0.886063	&	44.27					\\
42	&	50	&	50	&	100	&	50	&	0.656456	&	0.837083	&	1.073918	&	43.27					\\
43	&	50	&	50	&	100	&	50	&	0.75742	&	0.844813	&	0.963771	&	42.35					\\
44	&	50	&	50	&	150	&	50	&	0.67609	&	0.907669	&	1.180604	&	45.90					\\
45	&	50	&	50	&	150	&	50	&	0.671718	&	0.897607	&	1.117394	&	50.68					\\
46	&	50	&	50	&	200	&	50	&	0.664089	&	0.816423	&	1.043968	&	40.10					\\
47	&	50	&	50	&	200	&	50	&	0.571399	&	0.856775	&	1.074	&	56.78					\\
   \hline
48	&	50	&	75	&	50	&	50	&	0.514044	&	0.699661	&	0.899591	&	48.14					\\
49	&	50	&	75	&	50	&	50	&	0.647895	&	0.741642	&	0.978757	&	28.33					\\
50	&	50	&	75	&	100	&	50	&	0.623007	&	0.861603	&	1.058687	&	54.76					\\
51	&	50	&	75	&	100	&	50	&	0.607434	&	0.803321	&	0.940491	&	58.81					\\
52	&	50	&	75	&	150	&	50	&	0.742589	&	1.022707	&	1.157709	&	67.48					\\
53	&	50	&	75	&	150	&	50	&	0.706354	&	0.844787	&	0.980656	&	50.47					\\
54	&	50	&	75	&	200	&	50	&	0.719345	&	0.877183	&	1.064428	&	45.74					\\
55	&	50	&	75	&	200	&	50	&	0.690228	&	0.845069	&	0.968196	&	55.70					\\
  \hline
56	&	50	&	100	&	50	&	50	&	0.642142	&	0.838403	&	1.021128	&	51.79					\\
57	&	50	&	100	&	50	&	50	&	0.530841	&	0.864219	&	1.13489	&	55.19					\\
58	&	50	&	100	&	100	&	50	&	0.647984	&	0.832686	&	1.119443	&	39.18					\\
59	&	50	&	100	&	100	&	50	&	0.70828	&	0.925214	&	1.182871	&	45.71					\\
60	&	50	&	100	&	150	&	50	&	0.544611	&	0.755445	&	1.086004	&	38.94					\\
61	&	50	&	100	&	150	&	50	&	0.71695	&	0.891654	&	0.996318	&	62.54					\\
62	&	50	&	100	&	200	&	50	&	0.48787	&	0.771204	&	1.251992	&	37.08					\\
63	&	50	&	100	&	200	&	50	&	0.64742	&	0.812777	&	0.970732	&	51.14					\\
  \hline
\end{tabular}
}
\end{center}
\caption{Values of bilevel instances with dimension of $v$ equal to 50.
\label{table:vals_bil_50}}
\end{table}

\begin{table}
\begin{center}
{
\begin{tabular}{|r|r|r|r|r|r|r|r|r|r|r|r|r|}
\hline
& \multicolumn{3}{c|}{Dimensions} && LP & Preprocess & Optimal & \% Gap & \multicolumn{2}{c|}{Final \% gaps} \\ \cline{2-4} \cline{10-11}
\#	&	$v$	&	$b$	&	$g$	&	rank($Q$)	&	relaxation	&	lower bound	&	value	& shrunk	& CPLEX & Our code \\ \hline
64	&	75	&	25	&	50	&	50	&	0.70269	&	0.796622	&	1.015407	&	30.04	&	5.42	&		\\
65	&	75	&	25	&	50	&	50	&	0.523934	&	0.662829	&	0.925725	&	34.57	&	solved	&		\\
66	&	75	&	25	&	100	&	50	&	0.590782	&	0.727175	&	1.04984	&	29.71	&	17.98	&		\\
67	&	75	&	25	&	100	&	50	&	0.619884	&	0.743261	&	1.056703	&	28.24	&	18.16	&		\\  \hline
68	&	75	&	50	&	50	&	50	&	0.561663	&	0.78647	&	1.011143	&	50.01	&	solved	&		\\
69	&	75	&	50	&	50	&	50	&	0.560373	&	0.695674	&	0.980867	&	32.18	&	19.20	&		\\
70	&	75	&	50	&	100	&	50	&	0.801978	&	0.889381	&	1.041901	&	36.43	&	solved	&		\\
71	&	75	&	50	&	100	&	50	&	0.754207	&	0.881563	&	0.959749	&	61.96	&	solved	&		\\  \hline
72	&	75	&	75	&	50	&	50	&	0.700723	&	0.815533	&	1.026476	&	35.24	&	7.71	&		\\
73	&	75	&	75	&	50	&	50	&	0.601257	&	0.699462	&	0.930676	&	29.81	&	6.95	&		\\
74	&	75	&	75	&	100	&	50	&	0.346691	&	0.519392	&	0.86743	&	33.16	&	29.05	&		\\
75	&	75	&	75	&	100	&	50	&	0.594098	&	0.738111	&	1.094845	&	28.76	&	25.26	&		\\  \hline
76	&	100	&	25	&	50	&	50	&	0.882787	&	0.94306	&	0.992266	&	55.05	&	solved	&		\\
77	&	100	&	25	&	50	&	50	&	0.567423	&	0.613049	&	0.732117	&	27.70	&	3.08	&		\\
78	&	100	&	25	&	75	&	50	&	0.603912	&	0.639672	&	0.842811	&	14.97	&	20.18	&		\\
   \hline
79	&	75	&	25	&	50	&	75	&	0.475835	&	0.699877	&		&		&	no UB	&	5.19	\\
80	&	75	&	25	&	50	&	75	&	0.524895	&	0.732925	&	0.98902	&	44.82	&	7.86	&		\\
81	&	75	&	25	&	100	&	75	&	0.603074	&	0.783103	&		&		&	16.26	&	12.05	\\
82	&	75	&	25	&	100	&	75	&	0.512566	&	0.680474	&		&		&	16.44	&	12.09	\\
83	&	75	&	50	&	50	&	75	&	0.546256	&	0.697274	&	1.092166	&	27.66	&	25.04	&		\\
84	&	75	&	50	&	50	&	75	&	0.519369	&	0.68752	&		&		&	18.61	&	10.53	\\
85	&	75	&	50	&	100	&	75	&	0.508331	&	0.6817	&		&		&	22.22	&	5.29	\\
86	&	75	&	50	&	100	&	75	&	0.619981	&	0.80781	&		&		&	16.37	&	4.28	\\
87	&	75	&	75	&	50	&	75	&	0.485787	&	0.630083	&		&		&	27.69	&	13.95	\\
88	&	75	&	75	&	50	&	75	&	0.553843	&	0.737812	&		&		&	7.31	&	7.04	\\
89	&	75	&	75	&	100	&	75	&	0.416464	&	0.615691	&		&		&	18.33	&	13.76	\\
90	&	75	&	75	&	100	&	75	&	0.63856	&	0.783523	&	1.02771	&	37.25	&	solved	&		\\
   \hline
\end{tabular}
}
\end{center}
\caption{Values of bilevel instances with larger dimensions of~$v$.
The final gaps obtained by each code are indicated for
the instances it did not solve.
\label{table:vals_bil_larger}}
\end{table}

\begin{table}
\begin{center}
{\scriptsize  
\begin{tabular}{|r|r|r|r|r|r|r|r|r|r|}
\hline
& \multicolumn{3}{c|}{Dimensions} && Preprocess &  \multicolumn{2}{c|}{Our code} &
\multicolumn{2}{c|}{default CPLEX} \\ \cline{2-4} \cline{7-10}
\#	&	$v$	&	$b$	&	$g$	&	rank($Q$)	&	time	&	time	&	nodes	&	time	&	nodes 		\\ \hline
1       &       50      &       25      &       50      &       25      &       4.03    &       0.53    &       152     &       3.52    &       1395            \\
2       &       50      &       25      &       50      &       25      &       2.51    &       3.66    &       834     &       28.79   &       8785            \\
3       &       50      &       25      &       100     &       25      &       5.50    &       1.13    &       154     &       36.29   &       21962           \\
4       &       50      &       25      &       100     &       25      &       7.06    &       11.90   &       2706    &       364.99  &       512865          \\
5      &       50      &       25      &       150     &       25      &       10.10   &       4.31    &       974     &       39.37   &       20765           \\
6      &       50      &       25      &       150     &       25      &       14.56   &       16.56   &       4622    &       1072.49 &       966602          \\
7      &       50      &       25      &       200     &       25      &       19.37   &       9.09    &       2284    &       706.38  &       463141          \\
8      &       50      &       25      &       200     &       25      &       25.27   &       116.34  &       16582   &       444.74  &       245964          \\
\hline
\multicolumn{5}{r|}{\mbox{Means:}} & 11.05 & 20.44 & 3538 & 337.07 & 280184 \\
\hline
9      &       50      &       50      &       50      &       25      &       2.86    &       1.75    &       308     &       5.17    &       1897            \\
10      &       50      &       50      &       50      &       25      &       3.89    &       2.68    &       872     &       43.18   &       22949           \\
11      &       50      &       50      &       100     &       25      &       5.86    &       5.53    &       1118    &       22.96   &       12752           \\
12      &       50      &       50      &       100     &       25      &       7.83    &       1.61    &       222     &       39.19   &       26333           \\
13      &       50      &       50      &       150     &       25      &       11.11   &       6.04    &       1356    &       84.34   &       38741           \\
14      &       50      &       50      &       150     &       25      &       13.93   &       16.21   &       3194    &       216.39  &       136792          \\
15      &       50      &       50      &       200     &       25      &       7.85    &       0.01    &       2       &       0.46    &       0               \\
16      &       50      &       50      &       200     &       25      &       19.89   &       38.52   &       8810    &       1966.71 &       1297588         \\
\hline
\multicolumn{5}{r|}{\mbox{Means:}} & 9.15 & 9.04 & 1985 & 297.30 & 192131 \\
\hline
17      &       50      &       75      &       50      &       25      &       3.68    &       3.07    &       652     &       8.96    &       7615            \\
18      &       50      &       75      &       50      &       25      &       3.11    &       3.07    &       672     &       78.15   &       64029           \\
19      &       50      &       75      &       100     &       25      &       8.48    &       9.26    &       2748    &       217.66  &       192307          \\
20      &       50      &       75      &       100     &       25      &       8.81    &       12.40   &       3996    &       1035.67 &       833881          \\
21      &       50      &       75      &       150     &       25      &       5.71    &       0.04    &       14      &       0.43    &       0               \\
22      &       50      &       75      &       150     &       25      &       16.55   &       21.92   &       3628    &       3143.11 &       2896080         \\
23      &       50      &       75      &       200     &       25      &       21.39   &       17.50   &       2552    &       2861.97 &       1606015         \\
24      &       50      &       75      &       200     &       25      &       21.42   &       6.18    &       692     &       26.36   &       6211            \\ 
\hline
\multicolumn{5}{r|}{\mbox{Means:}} & 11.14 & 9.18 & 1869 & 921.54 & 700767 \\
\hline
25      &       50      &       100     &       50      &       25      &       4.52    &       8.85    &       1636    &       59.19   &       46175           \\
26      &       50      &       100     &       50      &       25      &       5.36    &       9.94    &       2614    &       234.31  &       224258          \\
27      &       50      &       100     &       100     &       25      &       5.81    &       5.98    &       716     &       49.56   &       31727           \\
28      &       50      &       100     &       150     &       25      &       13.10   &       27.25   &       4164    &       148.25  &       45614           \\
29      &       50      &       100     &       150     &       25      &       8.00    &       0.97    &       42      &       9.44    &       2528            \\
30      &       50      &       100     &       200     &       25      &       23.57   &       48.67   &       10332   &       3046.81 &       1566417         \\
31      &       50      &       100     &       200     &       25      &       13.62   &       8.13    &       824     &       299.12  &       135438          \\
\hline
\multicolumn{5}{r|}{\mbox{Means:}} & 10.57 & 15.68 & 2904 & 549.53 & 293165 \\
\hline
32       &       50      &       25      &       50      &       50      &       2.92    &       5.96    &       1556    &       64.01   &       36311           \\
33     &       50      &       25      &       50      &       50      &       2.70    &       7.72    &       2560    &       70.79   &       55449           \\
34     &       50      &       25      &       100     &       50      &       7.66    &       46.89   &       12238   &       112.51  &       56827           \\
35       &       50      &       25      &       100     &       50      &       8.15    &       29.99   &       5240    &       840.10  &       686855          \\
36       &       50      &       25      &       150     &       50      &       14.20   &       34.38   &       7064    &       $\geq$ 3600     &       2042811         \\
37      &       50      &       25      &       150     &       50      &       9.58    &       20.46   &       5170    &       255.14  &       74572           \\
38      &       50      &       25      &       200     &       50      &       17.34   &       161.83  &       34378   &       $\geq$ 3600     &       2001153         \\
39      &       50      &       25      &       200     &       50      &       18.56   &       168.98  &       39694   &       $\geq$ 3600     &       3371704         \\
\hline
\multicolumn{5}{r|}{\mbox{Means:}} & 10.14 & 59.53 & 13488  \\
  \hline
40      &       50      &       50      &       50      &       50      &       3.41    &       51.01   &       20124   &       283.70  &       277842          \\
41      &       50      &       50      &       50      &       50      &       3.61    &       9.07    &       2420    &       93.89   &       79268           \\
42      &       50      &       50      &       100     &       50      &       7.53    &       29.02   &       7926    &       $\geq$ 3600     &       4014451         \\
43      &       50      &       50      &       100     &       50      &       7.43    &       13.78   &       4422    &       208.67  &       71494           \\
44      &       50      &       50      &       150     &       50      &       13.48   &       96.14   &       14774   &       361.15  &       223345          \\
45      &       50      &       50      &       150     &       50      &       15.29   &       153.99  &       21876   &       581.26  &       360898          \\
46      &       50      &       50      &       200     &       50      &       17.46   &       248.31  &       48162   &       1008.91 &       465121          \\
47      &       50      &       50      &       200     &       50      &       25.74   &       143.38  &       26886   &       1165.71 &       648962          \\
\hline
\multicolumn{5}{r|}{\mbox{Means:}} & 11.74 & 93.09 & 18324  \\
    \hline
48      &       50      &       75      &       50      &       50      &       3.48    &       9.82    &       3334    &       43.88   &       23923           \\
49      &       50      &       75      &       50      &       50      &       3.67    &       22.58   &       6612    &       131.80  &       87870           \\
50      &       50      &       75      &       100     &       50      &       9.70    &       14.07   &       2166    &       28.00   &       8436            \\
51      &       50      &       75      &       100     &       50      &       7.12    &       4.84    &       662     &       101.97  &       46016           \\
52      &       50      &       75      &       150     &       50      &       13.43   &       2.39    &       232     &       3.83    &       926             \\
53      &       50      &       75      &       150     &       50      &       12.51   &       32.03   &       5386    &       178.99  &       118337          \\
54      &       50      &       75      &       200     &       50      &       21.96   &       56.51   &       5864    &       1881.67 &       904828          \\
55      &       50      &       75      &       200     &       50      &       21.35   &       26.14   &       4890    &       194.52  &       64333           \\
\hline
\multicolumn{5}{r|}{\mbox{Means:}} & 11.65 & 21.05 & 3643   \\
   \hline
56      &       50      &       100     &       50      &       50      &       4.87    &       13.45   &       3284    &       272.87  &       190537          \\
57      &       50      &       100     &       50      &       50      &       5.85    &       51.07   &       19014   &       $\geq$ 3600     &       30944839                \\
58      &       50      &       100     &       100     &       50      &       8.08    &       120.59  &       33178   &       2334.10 &       1597429         \\
59      &       50      &       100     &       100     &       50      &       9.86    &       147.33  &       28058   &       2008.30 &       1283510         \\
60      &       50      &       100     &       150     &       50      &       16.61   &       218.98  &       48542   &       $\geq$ 3600     &       3310389         \\
61      &       50      &       100     &       150     &       50      &       11.40   &       18.82   &       2662    &       46.03   &       13536           \\
62      &       50      &       100     &       200     &       50      &       21.63   &       3122.24 &       323810  &       $\geq$ 3600     &       13525363                \\
63      &       50      &       100     &       200     &       50      &       22.05   &       91.29   &       9660    &       353.94  &       122675          \\
\hline
\multicolumn{5}{r|}{\mbox{Means:}} & 12.54 & 472.97 & 58526   \\
\cline{6-8}
\end{tabular}
}
\end{center}
\caption{Performance on bilevel instances with dimension of $v$ equal to 50.
\label{table:perf_bil_50}}
\end{table}

\begin{table}
\begin{center}
{
\begin{tabular}{|r|r|r|r|r|r|r|r|r|r|}
\hline
& \multicolumn{3}{c|}{Dimensions} && Preprocess &  \multicolumn{2}{c|}{Our code} &
\multicolumn{2}{c|}{default CPLEX} \\ \cline{2-4} \cline{7-10}
\#	&	$v$	&	$b$	&	$g$	&	rank($Q$)	&	time	&	time	&	nodes	&	time	&	nodes 		\\ \hline
64      &       75      &       25      &       50      &       50      &       11.47   &       290.70  &       75616   &       $\geq$ 3600     &       12455384                \\
65      &       75      &       25      &       50      &       50      &       7.73    &       299.37  &       68770   &       $\geq$ 3600     &       17447435                \\
66      &       75      &       25      &       100     &       50      &       20.20   &       396.67  &       86638   &       $\geq$ 3600     &       14159416                \\
67      &       75      &       25      &       100     &       50      &       16.82   &       934.03  &       159358  &       $\geq$ 3600     &       7263261         \\
\hline
\multicolumn{5}{r|}{\mbox{Means:}} & 14.06 & 480.19 & 97596 \\
\hline
68      &       75      &       50      &       50      &       50      &       12.75   &       89.09   &       19502   &       $\geq$ 3600     &       31475242                \\
69      &       75      &       50      &       50      &       50      &       10.46   &       139.84  &       28866   &       $\geq$ 3600     &       1747600         \\
70      &       75      &       50      &       100     &       50      &       13.46   &       51.33   &       7478    &       $\geq$ 3600     &       4579441         \\
71      &       75      &       50      &       100     &       50      &       17.91   &       5.75    &       362     &       3585.51 &       1415201         \\
\hline
\multicolumn{5}{r|}{\mbox{Means:}} & 13.65 & 71.50 & 14052 \\
\hline
72      &       75      &       75      &       50      &       50      &       9.32    &       94.15   &       16830   &       $\geq$ 3600     &       4006474         \\
73      &       75      &       75      &       50      &       50      &       11.50   &       86.90   &       15076   &       $\geq$ 3600     &       3636937         \\
74      &       75      &       75      &       100     &       50      &       21.87   &       2519.64 &       267744  &       $\geq$ 3600     &       1185965         \\
75      &       75      &       75      &       100     &       50      &       24.30   &       1569.69 &       197996  &       $\geq$ 3600     &       1144374         \\
\hline
\multicolumn{5}{r|}{\mbox{Means:}} & 16.75 & 1067.60 & 124412 \\
\hline
76      &       100     &       25      &       50      &       50      &       8.99    &       6.78    &       408     &       155.29  &       27742           \\
77      &       100     &       25      &       50      &       50      &       8.18    &       79.48   &       5450    &       $\geq$ 3600     &       2647874         \\ 
78      &       100     &       25      &       75      &       50      &       26.23   &       1030.75 &       100044  &       $\geq$ 3600     &       1049577         \\
\hline
\multicolumn{5}{r|}{\mbox{Means:}} & 14.47 & 372.34 & 35301 \\
\hline
79      &       75      &       25      &       50      &       75      &       8.50    &       $\geq$ 3600     &       542572  &       $\geq$ 3600     &       3248691         \\
80      &       75      &       25      &       50      &       75      &       11.67   &       1197.30 &       212828  &       $\geq$ 3600     &       2730181         \\
81      &       75      &       25      &       100     &       75      &       16.93   &       $\geq$ 3600     &       449218  &       $\geq$ 3600     &       13374465                \\
82      &       75      &       25      &       100     &       75      &       22.19   &       $\geq$ 3600     &       268486  &       $\geq$ 3600     &       13553753                \\
83      &       75      &       50      &       50      &       75      &       10.94   &       3587.62 &       445828  &       $\geq$ 3600     &       2126574         \\
84      &       75      &       50      &       50      &       75      &       12.60   &       $\geq$ 3600     &       488384  &       $\geq$ 3600     &       1758366         \\
85      &       75      &       50      &       100     &       75      &       22.05   &       $\geq$ 3600     &       408186  &       $\geq$ 3600     &       1269948         \\
86      &       75      &       50      &       100     &       75      &       19.24   &       $\geq$ 3600     &       429360  &       $\geq$ 3600     &       1389433         \\
87      &       75      &       75      &       50      &       75      &       12.25   &       $\geq$ 3600     &       442694  &       $\geq$ 3600     &       1984975         \\
88      &       75      &       75      &       50      &       75      &       11.15   &       $\geq$ 3600     &       463366  &       $\geq$ 3600     &       2227245         \\
89      &       75      &       75      &       100     &       75      &       24.20   &       $\geq$ 3600     &       276876  &       $\geq$ 3600     &       1390365         \\
90      &       75      &       75      &       100     &       75      &       20.69   &       1221.35 &       107572  &       3537.29 &       1604647         \\
\hline
\end{tabular}
}
\end{center}
\caption{Performance on bilevel instances with larger dimensions of~$v$.
\label{table:perf_bil_larger}}
\end{table}

\section{Inverse QP Instances  \label{app:inverse}}

Computational results on 25 inverse QP instances can be found in Table~\ref{table.inverseQP}.
For each set of 5 instances, the average CPU time is listed if all the instances were solved
or the number of solved instances is noted.

\begin{table}
\begin{center}
{
\begin{tabular}{|r|r|r|r|r|r|r|r|r|}
\hline
$\tm$ & $\tn$ & instance & time for our code & CPLEX MIP time
& \multicolumn{4}{c|}{CPLEX MIP cuts} \\ \cline{6-9}
&&&&& GF & MIR & L\&P & IB \\ \hline \hline
100	&	75	&	a	&	303.58	&	3504.84	&5&3&2&\\
100	&	75	&	b	&	304.58	&	2686.84	&22&1&&\\
100	&	75	&	c	&	60.38	&	69.09	&11&&&1\\
100	&	75	&	d	&	328.71	&	3600.00	&19&&&\\
100	&	75	&	e	&	94.91	&	182.61	&11&1&&\\ \hline
\multicolumn{3}{r|}{mean or success}&	218.43	&	4 of 5	\\ \hline \hline
120	&	90	&	a	&	271.24	&	1314.80	&9&7&&\\
120	&	90	&	b	&	250.26	&	3600.00	&14&2&&\\
120	&	90	&	c	&	126.13	&	209.95	&7&5&&\\
120	&	90	&	d	&	215.09	&	3600.00	&16&&&\\
120	&	90	&	e	&	1101.49	&	3864.48	&11&8&&\\ \hline
\multicolumn{3}{r|}{mean or success}&	392.84	&	2 of 5	\\ \hline \hline
150	&	20	&	a	&	55.12	&	702.95	&5&&&\\
150	&	20	&	b	&	486.00	&	3263.58	&7&&&\\
150	&	20	&	c	&	163.54	&	3319.01	&4&&&\\
150	&	20	&	d	&	49.13	&	1260.59	&6&&&\\
150	&	20	&	e	&	3605.77	&	3781.06	&1&&&\\ \hline
\multicolumn{3}{r|}{success}		&	4 of 5	&	4 of 5	\\ \hline \hline
200	&	15	&	a	&	154.11	&	1057.09	&7&&&\\
200	&	15	&	b	&	81.65	&	477.30	&1&&&\\
200	&	15	&	c	&	3604.85	&	3600.00	&4&&&\\
200	&	15	&	d	&	128.41	&	365.59	&7&&1&\\
200	&	15	&	e	&	47.71	&	224.61	&5&&&\\ \hline
\multicolumn{3}{r|}{success}		&	4 of 5	&	4 of 5	\\ \hline \hline
400	&	5	&	a	&	225.54	&	1573.05	&7&&&\\
400	&	5	&	b	&	1797.89	&	3600.00	&17&&&\\
400	&	5	&	c	&	238.79	&	183.92	&21&&&\\
400	&	5	&	d	&	252.45	&	388.98	&20&&&\\
400	&	5	&	e	&	47.45	&	224.18	&5&&&\\ \hline
\multicolumn{3}{r|}{mean or success}&	512.42	&	4 of 5	\\ \cline{4-5}
\end{tabular}
}
\end{center}
\caption{Performance on 25 inverse quadratic programs.
Mean solution time is listed for each set of five problems solved successfully by a code;
otherwise, the number of solved instances is given.
The number of cutting planes added by CPLEX MIP is also reported;
GF are Gomory fractional cuts, MIR are mixed integer rounding cuts,
L\&P are lift-and-project cuts, and IB are implicit bound cuts.
\label{table.inverseQP}}
\end{table}

\end{document}